\documentclass{amsart}

\usepackage{amsfonts}
\usepackage{amssymb}
\usepackage{graphicx}
\usepackage{stmaryrd}
\usepackage{psfrag}
\usepackage{xypic}

\newcommand{\N}{{\mathbb N}}
\newcommand{\Z}{{\mathbb Z}}
\newcommand{\Q}{{\mathbb Q}}

\newcommand{\C}{{\mathbb C}}

\newcommand{\RR}{{\mathcal R}}
\newcommand{\TT}{{\mathcal T}}
\newcommand{\CC}{{\mathcal C}}
\newcommand{\LL}{{\mathcal L}}
\newcommand{\GG}{{\mathcal G}}

\newcommand{\Spec}{\operatorname{Spec}}
\newcommand{\Schub}{{\mathfrak{S}}}
\newcommand{\Groth}{{\mathfrak{G}}}
\newcommand{\bull}{{\sssize \bullet}}

\newcommand{\rank}{\operatorname{rank}}
\newcommand{\sh}{\operatorname{sh}}

\newcommand{\Gr}{\operatorname{Gr}}
\newcommand{\Fl}{\operatorname{F\ell}}
\renewcommand{\O}{{\mathcal O}}
\newcommand{\I}{{\mathcal I}}
\newcommand{\col}{\operatorname{col}}
\newcommand{\GL}{\operatorname{GL}}

\newcommand{\pic}[2]{\includegraphics[scale=0.#1]{#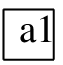}}

\newcommand{\tableau}[1]{\includegraphics{#1.eps}}
\newcommand{\rtab}[2]{\raisebox{#1}{\includegraphics{#2.eps}}}

\newcommand{\ins}[1]{\xrightarrow[#1]{}}

\newtheorem{thm}{Theorem}[section] %\numberwithin{thm}{chapter}
\newtheorem{lemma}[thm]{Lemma} 
\newtheorem{prop}[thm]{Proposition}
\newtheorem{cor}[thm]{Corollary}
\newtheorem{conj}[thm]{Conjecture}

\theoremstyle{definition}
\newtheorem{remark}[thm]{Remark}
\newtheorem{example}[thm]{Example}
\newtheorem{defn}[thm]{Definition}

\newcommand{\refsec}[1]{Section~\ref{#1}}
\newcommand{\refeqn}[1]{(\ref{#1})}
\newcommand{\refthm}[1]{Theorem~\ref{#1}}
\newcommand{\refprop}[1]{Proposition~\ref{#1}}
\newcommand{\reflemma}[1]{Lemma~\ref{#1}}
\newcommand{\refcor}[1]{Corollary~\ref{#1}}
\newcommand{\refconj}[1]{Conjecture~\ref{#1}}
\newcommand{\refdefn}[1]{Definition~\ref{#1}}
\newcommand{\refexm}[1]{Example~\ref{#1}}

\newenvironment{romenum}{\begin{enumerate}}{\end{enumerate}}

\psfrag{tht1}{$\theta_1$}
\psfrag{tht2}{$\theta_2$}
\psfrag{lambda}{$\lambda$}
\psfrag{sigma}{$\sigma$}
\psfrag{mu}{$\mu$}
\psfrag{t1}{$\theta_1$}
\psfrag{t2}{$\theta_2$}
\psfrag{t3}{$\theta_3$}
\psfrag{t4}{$\theta_4$}
\psfrag{t5}{$\theta_5$}

\psfrag{ac1}{$123$}
\psfrag{aa1}{$1$}
\psfrag{aa2}{$2\,3$}
\psfrag{aa3}{$2$}
\psfrag{ad1}{$1$}
\psfrag{ad2}{$2$}
\psfrag{ad3}{$2$}
\psfrag{ad4}{$2$}
\psfrag{ad5}{$3$}
\psfrag{b1}{$1$}
\psfrag{b2}{$2$}
\psfrag{a1}{$2$}
\psfrag{ag1}{$1$}
\psfrag{ag2}{$2$}
\psfrag{ag3}{$2$}
\psfrag{ag4}{$2\,3$}
\psfrag{ah1}{$2\,3$}
\psfrag{ai1}{$1\,2$}
\psfrag{ai2}{$2$}
\psfrag{ai3}{$3$}
\psfrag{aj1}{$2$}
\psfrag{aj2}{$2$}
\psfrag{aj3}{$3$}
\psfrag{bb1}{$1\,2$}
\psfrag{bb2}{$4\,5$}
\psfrag{bc1}{$1$}
\psfrag{bc2}{$245$}
\psfrag{bc3}{$2\,3$}
\psfrag{bc4}{$5$}
\psfrag{ak1}{$a$}
\psfrag{al1}{$a$}
\psfrag{al2}{$x$}
\psfrag{am1}{$a\,b$}
\psfrag{an1}{$b$}
\psfrag{an2}{$a$}
\psfrag{an3}{$x$}
\psfrag{ao1}{$a\,b$}
\psfrag{ao2}{$c$}
\psfrag{ap1}{$b$}
\psfrag{ap2}{$a$}
\psfrag{ap3}{$x\,c$}
\psfrag{aq1}{$a$}
\psfrag{aq2}{$b$}
\psfrag{ar1}{$b$}
\psfrag{ar2}{$a$}
\psfrag{ar3}{$x$}
\psfrag{as1}{$b$}
\psfrag{at1}{$x$}
\psfrag{at2}{$b$}
\psfrag{av1}{$a$}
\psfrag{av2}{$x\,b$}
\psfrag{au1}{$x\,b$}
\psfrag{da1}{$a$}
\psfrag{re1}{$y$}
\psfrag{re2}{$a$}
\psfrag{re3}{$b\,c$}
\psfrag{rf1}{$a\,y$}
\psfrag{rf2}{$c$}
\psfrag{ra1}{$y$}
\psfrag{ra2}{$a$}
\psfrag{ra3}{$b$}
\psfrag{rb1}{$y$}
\psfrag{rb2}{$b$}
\psfrag{rc1}{$y$}
\psfrag{rc2}{$a$}
\psfrag{rd1}{$y$}
\psfrag{rg1}{$a\,y$}
\psfrag{rg2}{$b$}
\psfrag{rh1}{$a\,y$}
\psfrag{c1}{$a\,b$}
\psfrag{az1}{$C_1$}
\psfrag{az2}{$C_2$}
\psfrag{az3}{$\cdots$}
\psfrag{az4}{$C_\ell$}
\psfrag{bm1}{$1$}
\psfrag{bm2}{$2\,3$}
\psfrag{bm3}{$1\,2$}
\psfrag{bm4}{$2 3 4$}
\psfrag{bm5}{$2$}
\psfrag{bm6}{$3\,5$}
\psfrag{bm7}{$7$}
\psfrag{bn1}{1}
\psfrag{bn2}{1}
\psfrag{bo1}{1}
\psfrag{bo2}{2}
\psfrag{bp1}{1}
\psfrag{bp2}{$1\,2$}
\psfrag{bh1}{$1$}
\psfrag{bh2}{$1$}
\psfrag{bh3}{$2$}
\psfrag{bh4}{$2$}
\psfrag{bg1}{$1$}
\psfrag{bg2}{$1$}
\psfrag{bg3}{$3$}
\psfrag{bg4}{$2$}
\psfrag{bg5}{$2$}
\psfrag{bg6}{$4$}
\psfrag{bg7}{$3$}
\psfrag{bg8}{$4$}
\psfrag{bg9}{$5$}
\psfrag{bd1}{$1\,2$}
\psfrag{bd2}{$2$}
\psfrag{bd3}{$2\,3$}
\psfrag{bd4}{$4$}
\psfrag{be1}{$1$}
\psfrag{be2}{$2$}
\psfrag{be3}{$2$}
\psfrag{be4}{$4$}
\psfrag{be5}{$2$}
\psfrag{be6}{$3$}
\psfrag{ca1}{$a$}
\psfrag{ca2}{$b$}
\psfrag{dz1}{6}
\psfrag{dz2}{7}
\psfrag{dz3}{8}
\psfrag{dz4}{4}
\psfrag{dz5}{5}
\psfrag{dz6}{6}
\psfrag{dz7}{3}
\psfrag{dz8}{4}
\psfrag{bj1}{$3\,5$}
\psfrag{bj2}{$5\,6$}
\psfrag{bj3}{$2$}
\psfrag{bj4}{$6$}
\psfrag{bj5}{$7$}
\psfrag{bj6}{$2\,3$}
\psfrag{bj7}{$3$}
\psfrag{bj8}{$7$}
\psfrag{bj9}{$9$}
\psfrag{bj10}{$4$}
\psfrag{bj11}{$4\,5$}
\psfrag{bj12}{$8$}
\psfrag{bj13}{$8$}
\psfrag{bj14}{$8$}
\psfrag{bk1}{$3\,5$}
\psfrag{bk2}{$5$}
\psfrag{bk3}{$2$}
\psfrag{bk4}{$6$}
\psfrag{bk5}{$6$}
\psfrag{bk6}{$2\,3$}
\psfrag{bk7}{$3$}
\psfrag{bk8}{$7$}
\psfrag{bk9}{$9$}
\psfrag{bk10}{$4$}
\psfrag{bk11}{$4\,5$}
\psfrag{bk12}{$8$}
\psfrag{bk13}{$7$}
\psfrag{bk14}{$8$}
\psfrag{bf1}{$1$}
\psfrag{bf2}{$1\,2$}
\psfrag{bf3}{$3\,4$}
\psfrag{bf4}{$2\,3$}
\psfrag{bf5}{$4\,5$}
\psfrag{dc1}{$y_k$}
\psfrag{dc2}{$C_k$}
\psfrag{dd1}{$y$}
\psfrag{dd2}{$C_1$}
\psfrag{de1}{$y_\ell$}
\psfrag{de2}{$C_1'$}
\psfrag{de3}{$C_2'$}
\psfrag{de4}{$\cdots$}
\psfrag{de5}{$C_\ell'$}
\psfrag{ba1}{$235$}
\psfrag{ay1}{$\mu$}
\psfrag{ay2}{$\lambda$}
\psfrag{bi1}{$\rho$}
\psfrag{bi2}{$\lambda$}
\psfrag{bi3}{$\mu$}
\psfrag{bl1}{$T_1$}
\psfrag{bl2}{$D$}
\psfrag{bl3}{$T_2$}
\psfrag{bl4}{$a$}
\psfrag{bl5}{$C$}

\numberwithin{equation}{section}

\begin{document}

\title[A Littlewood-Richardson rule for the $K$-theory of Grassmannians]{
  A Littlewood-Richardson rule for the \\ 
  $K$-theory of Grassmannians} 
\author{Anders Skovsted Buch}
%\date{\today}
\address{Massachusetts Institute of Technology \\
  Building 2, Room 275 \\
  77 Massachusetts Avenue \\
  Cambridge, MA 02139
}
\email{abuch@math.mit.edu}
\thanks{The author was partially supported by NSF Grant DMS-0070479}
\thanks{{\em E-mail address}: abuch@math.mit.edu}

\begin{abstract}
  We prove an explicit combinatorial formula for the structure
  constants of the Grothendieck ring of a Grassmann variety with
  respect to its basis of Schubert structure sheaves.  We furthermore
  relate $K$-theory of Grassmannians to a bialgebra of stable
  Grothendieck polynomials, which is a $K$-theory parallel of the ring
  of symmetric functions.
\end{abstract}

\maketitle

\section{Introduction}

Let $X = \Gr(d,\C^n)$ be the Grassmann variety of $d$ dimensional
subspaces of $\C^n$.  The goal of this paper is to give an explicit
combinatorial description of the Grothendieck ring $K^\circ X$ of
algebraic vector bundles on $X$.

$K$-theory of Grassmannians is a special case of $K$-theory of flag
varieties, which was studied by Kostant and Kumar
\cite{kostant.kumar:t-equivariant} and by Demazure
\cite{demazure:desingularisation}.  Lascoux and Sch{\"u}tzenberger
defined Grothendieck polynomials which give formulas for the structure
sheaves of the Schubert varieties in a flag variety
\cite{lascoux.schutzenberger:structure, lascoux:anneau}.  The
combinatorial understanding of these polynomials was further developed
by Fomin and Kirillov \cite{fomin.kirillov:yang-baxter,
  fomin.kirillov:grothendieck}.

Recall that if $\lambda = (\lambda_1 \geq \lambda_2 \geq \cdots \geq
\lambda_d)$ is a partition with $d$ parts and $\lambda_1 \leq n-d$,
then the Schubert variety in $X$ associated to $\lambda$ is the subset
\begin{equation}
\label{eqn_grschub}
  \Omega_\lambda = \{ V \in \Gr(d,\C^n) \mid 
  \dim(V \cap \C^{n-d+i-\lambda_i}) \geq i ~\forall 1 \leq i \leq d \}
  \,.
\end{equation}
Here $\C^k \subset \C^n$ denotes the subset of vectors whose last
$n-k$ components are zero.  The codimension of $\Omega_\lambda$ is
equal to the weight $|\lambda| = \sum \lambda_i$ of $\lambda$.  If we
identify partitions with their Young diagrams, then a Schubert variety
$\Omega_\mu$ is contained in $\Omega_\lambda$ if and only $\mu$
contains $\lambda$.  From the fact that the open Schubert cells
$\Omega_\lambda^\circ = \Omega_\lambda \smallsetminus \cup_{\mu
  \varsupsetneq \lambda} \Omega_\mu$ form a cell decomposition of $X$,
one can deduce that the classes of the structure sheaves
$\O_{\Omega_\lambda}$ form a basis for the Grothendieck ring of $X$.

We will study the structure constants for $K^\circ X$ with respect to
this basis.  These are the unique integers $c^\nu_{\lambda \mu}$ such
that
\[ [\O_{\Omega_\lambda}] \cdot [\O_{\Omega_\mu}] = 
   \sum_\nu c^\nu_{\lambda \mu} \, [\O_{\Omega_\nu}] \,.
\]
These constants depend only on the partitions $\lambda$, $\mu$, and
$\nu$, not on the Grassmannian where the Schubert varieties for these
partitions are realized.  Furthermore $c^\nu_{\lambda \mu}$ is
non-zero only if $|\nu| \geq |\lambda| + |\mu|$.  The constants
$c^\nu_{\lambda \mu}$ for which $|\nu| = |\lambda| + |\mu|$ are the
usual Littlewood-Richardson coefficients, i.e.\ the structure
constants for the cohomology ring $H^*(X)$ with respect to the basis
of cohomology classes of Schubert varieties.  Another known case is a
Pieri formula of Lenart which expresses the coefficients
$c^\nu_{\lambda,(k)}$ for multiplying with the structure sheaf of a
special Schubert variety $\Omega_{(k)}$ as binomial coefficients
\cite{lenart:combinatorial}.  Notice that since the duality
isomorphism $\Gr(d,\C^n) \longrightarrow \Gr(n-d,\C^n)$ takes
$\Omega_\lambda$ to the Schubert variety $\Omega_{\lambda'}$ for the
conjugate partition \cite[Ex.~9.20]{fulton:young}, the structure
constants must satisfy $c^\nu_{\lambda \mu} = c^{\nu'}_{\lambda'
  \mu'}$.

Our main result is an explicit combinatorial formula stating that the
coefficient $c^\nu_{\lambda \mu}$ is $(-1)^{|\nu|-|\lambda|-|\mu|}$
times the number of objects called {\em set-valued tableaux\/} which
satisfy certain properties.  Set-valued tableaux are similar to
semistandard Young tableaux, but allow a non-empty set of integers in
each box of a Young diagram rather than a single integer.  When $|\nu|
= |\lambda| + |\mu|$ our formula specializes to the classical
Littlewood-Richardson rule.

Our formula implies that if $c^\nu_{\lambda \mu}$ is not zero, then
$\nu$ is contained in the union of all partitions $\rho$ of weight
$|\rho| = |\lambda| + |\mu|$ such that the Littlewood-Richardson
coefficient $c^\rho_{\lambda \mu}$ is non-zero.  As a consequence, for
fixed $\lambda$ and $\mu$ the constants $c^\nu_{\lambda \mu}$ are only
non-zero for finitely many partitions $\nu$, even though one might
conceivably get arbitrarily many such constants by realizing the
product $[\O_{\Omega_\lambda}] \cdot [\O_{\Omega_\mu}]$ in larger and
larger Grassmannians.  We do not know any geometric reason for this or
for the alternating signs of the structure constants.

This observation allows us to define a commutative ring $\Gamma =
\bigoplus \Z \cdot G_\lambda$ with a formal basis $\{G_\lambda\}$
indexed by partitions and multiplication defined by $G_\lambda \cdot
G_\mu = \sum_\nu c^\nu_{\lambda \mu} G_\nu$.  The Grothendieck ring
$K^\circ \Gr(d,\C^n)$ is then the quotient of this ring by the ideal
spanned by the basis elements $G_\lambda$ for partitions that do not
fit in a rectangle with $d$ rows and $n-d$ columns.  The pullbacks of
Grothendieck rings defined by the natural embeddings $\Gr(d_1,
\C^{n_1}) \times \Gr(d_2, \C^{n_2}) \subset \Gr(d_1+d_2,
\C^{n_1+n_2})$ furthermore define a coproduct on $\Gamma$ which makes
it a bialgebra.

This bialgebra $\Gamma$ can be seen as a $K$-theory parallel of the
ring of symmetric functions \cite{macdonald:symmetric*2,
  fulton:young}, which in a similar way describes the cohomology of
Grassmannians, in addition to representation theory of symmetric and
general linear groups and numerous other areas.  Furthermore, if we
define a filtration of $\Gamma$ by ideals $\Gamma_p =
\bigoplus_{|\lambda| \geq p} \Z \cdot G_\lambda$, then the associated
graded bialgebra is naturally the ring of symmetric functions.  This
filtration corresponds to Grothendieck's $\gamma$-filtration of the
$K$-ring of any non-singular algebraic variety.  In general the
associated graded ring is isomorphic to the Chow ring of the variety
after tensoring with $\Q$.

We will realize the algebra $\Gamma$ as the linear span of all stable
Grothendieck polynomials, which we show has a basis indexed by
Grassmannian permutations.  The Littlewood-Richardson rule is proved
by defining stable Grothendieck polynomials in non-commutative
variables, and showing that these polynomials multiply exactly like
those in commutative variables.  In order to carry out this
construction we will define a {\em jeu de taquin\/} algorithm for
set-valued tableaux.

In \refsec{sec_gpoly} we fix the notation concerning Grothendieck
polynomials.  In \refsec{sec_settab} we then prove a formula for
stable Grothendieck polynomials of 321-avoiding permutations in terms
of set-valued tableaux.  This formula uses the skew diagram associated
to a 321-avoiding permutation \cite{billey.jockusch.ea:some} and is
derived from a more general formula for Grothendieck polynomials of
Fomin and Kirillov \cite{fomin.kirillov:yang-baxter}.
\refsec{sec_jdt} develops a column bumping algorithm for set-valued
tableaux, which in \refsec{sec_locplac} is used to prove the
Littlewood-Richardson rule for the structure constants $c^\nu_{\lambda
  \mu}$ in $\Gamma$.  In \refsec{sec_coprod} we derive similar
Littlewood-Richardson rules for the coproduct in $\Gamma$ and for
writing the stable Grothendieck polynomial of any 321-avoiding
permutation as a linear combination of the basis elements of $\Gamma$.
In \refsec{sec_conseq} we deduce a number of consequences of these
rules, including a Pieri formula for the coproduct in $\Gamma$, a
result about multiplicity free products, and the above described bound
on partitions $\nu$ for which $c^\nu_{\lambda \mu}$ is not zero.  In
\refsec{sec_geometry} the relationship between $\Gamma$ and $K$-theory
of Grassmannians is established.  In addition we use the methods
developed in this paper to give simple proofs of some unpublished
results of A.~Knutson regarding triple intersections in $K$-theory.
\refsec{sec_structure} finally contains a discussion of the overall
structure of the bialgebra $\Gamma$.  We show that if the inverse of
the element $t = 1 - G_1$ is joined to $\Gamma$ then the result
$\Gamma_t$ is a Hopf algebra.  We furthermore pose a conjecture which
implies that the Abelian group scheme $\Spec \Gamma_t$ looks like an
infinite affine space minus a hyperplane.  We conclude by raising some
additional questions.  We hope that the statements in the last two
sections will be comprehensible after reading \refsec{sec_gpoly} and
the first seven lines of \refsec{sec_coprod}.

This paper came out of a project aimed at finding a formula for the
structure sheaf of a quiver variety.  We will present such a formula
in \cite{buch:structure}, thus generalizing our earlier results with
W.~Fulton regarding the cohomology class of a quiver variety
\cite{buch.fulton:chern}.  We thank Fulton for numerous helpful
discussions and suggestions during the project.  For carrying out the
work in this paper, it has been invaluable to speak to S.~Fomin, from
whom we have learned very much about Grothendieck polynomials during
several fruitful discussions.  In particular we thank Fomin for
supplying the proof of \reflemma{lemma_conjugate} and for suggestions
that led to simplifications of many proofs.  We also thank F.~Sottile
for useful discussions about stable Grothendieck polynomials.  We are
grateful to A.~Knutson for sharing his ideas about $K$-theoretic
triple intersections and for allowing us to report about them here.
Finally, we thank A.~Lascoux for informing us about a remarkable
recursive formula for stable Grothendieck polynomials which implies
that our \refconj{conj_alpha} is true.

%%% Local Variables: 
%%% mode: latex
%%% TeX-master: "gamma"
%%% End: 

\section{Grothendieck polynomials}
\label{sec_gpoly}

% \subsection{Definitions}

In this section we fix the notation regarding Grothendieck polynomials
and stable Grothendieck polynomials.  Grothendieck polynomials were
introduced by Lascoux and Sch{\"u}tzenberger as representatives for
the structure sheaves of the Schubert varieties in a Flag variety
\cite{lascoux.schutzenberger:structure, lascoux:anneau}.  For any
permutation $w \in S_n$ we define the {\em double Grothendieck
  polynomial\/} $\Groth_w = \Groth_w(x;y)$ as follows.  If $w$ is the
longest permutation $w_0 = n \, (n-1) \cdots 2 \, 1$ we set
\[ \Groth_{w_0} = \prod_{i+j \leq n} (x_i + y_j - x_i \, y_j) \,. \]
Otherwise we can find a simple reflection $s_i = (i,i+1) \in S_n$ such
that $\ell(w \, s_i) = \ell(w) + 1$.  Here $\ell(w)$ denotes the
length of $w$, which is the smallest number $\ell$ for which $w$ can
be written as a product of $\ell$ simple reflections.  We then define
\[ \Groth_w = \pi_i(\Groth_{w s_i}) \]
where $\pi_i$ is the {\em isobaric divided difference operator\/}
given by
\[ \pi_i(f) = \frac{(1-x_{i+1}) f(x_1,x_2,\dots) - 
   (1-x_i) f(\dots,x_{i+1},x_i,\dots)}{x_i - x_{i+1}} \,.
\]
This definition is independent of our choice of the simple reflection
$s_i$ since the operators $\pi_i$ satisfy the Coxeter relations.

Notice that the longest element in $S_{n+1}$ is $w_0^{(n+1)} = w_0
\cdot s_n \cdot s_{n-1} \cdots s_1$.  Since $\pi_n \cdot \pi_{n-1}
\cdots \pi_1$ applied to the Grothendieck polynomial for $w_0^{(n+1)}$
is equal to $\Groth_{w_0}$, it follows that $\Groth_w$ does not depend
on which symmetric group $w$ is considered an element of.

Let $1^m \times w \in S_{m+n}$ denote the permutation which is the
identity on $\{1, 2, \dots, m\}$ and which maps $j$ to $w(j-m) + m$
for $j > m$.  Fomin and Kirillov have shown that when $m$ grows to
infinity, the coefficient of each fixed monomial in $\Groth_{1^m
  \times w}$ eventually becomes stable
\cite{fomin.kirillov:yang-baxter}.  The {\em stable Grothendieck
  polynomial\/} $G_w \in \Z\llbracket x_i,y_i \rrbracket_{i \geq 0}$
is defined as the resulting power series:
\[ G_w = G_w(x;y) = \lim_{m \to \infty} \Groth_{1^m \times w} \,. \]
Fomin and Kirillov also proved that this power series is symmetric in
the variables $\{x_i\}$ and $\{y_i\}$ separately, and that
\[ G_w(1-e^{-x}; 1-e^y) = 
   G_w(1-e^{-x_1}, 1-e^{-x_2}, \dots; 1-e^{y_1}, 1-e^{y_2}, \dots) 
\]
is super symmetric, i.e.\ if one sets $x_1 = y_1$ in this expression
then the result is independent of $x_1$ and $y_1$.

If we put all the variables $y_i$ equal to zero in $\Groth_w(x;y)$, we
obtain the single Grothendieck polynomial $\Groth_w(x) =
\Groth_w(x;0)$.  Similarly the single stable Grothen\-dieck polynomial
for $w$ is defined as $G_w(x) = G_w(x;0)$.  Notice that the super
symmetry of $G_w(1-e^{-x}; 1-e^y)$ implies that the double stable
Grothendieck polynomial $G_w(x;y)$ is uniquely determined by the
single polynomial $G_w(x)$ \cite{stembridge:characterization,
  macdonald:symmetric*2}.  We will use the notation $G_w(x_1, x_2,
\dots, x_p; y_1, y_2, \dots, y_q)$ for the polynomial obtained by
setting $x_i = 0$ for $i > p$ and $y_j = 0$ for $j > q$ in the stable
polynomial $G_w(x;y)$.

% \subsection{321-Avoiding permutations}

If $\lambda \subset \nu$ are two partitions, let $\nu/\lambda$ denote
the skew diagram of boxes in $\nu$ which are not in $\lambda$, and let
$|\nu/\lambda|$ be the number of boxes in this diagram.  Now choose
a numbering of the north-west to south-east diagonals in the diagram
with positive integers, which increase consecutively from south-west
to north-east.  For example, if $\nu = (4,3,2)$ and $\lambda = (1)$,
and if the bottom-left box in $\nu/\lambda$ is in diagonal number $3$,
then the numbering is given by the picture:
\[ \tableau{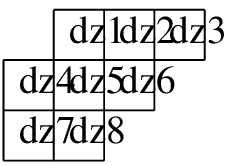} \]
Let $(i_1, i_2, \dots, i_m)$ be the sequence of diagonal numbers of
the boxes in $\nu/\lambda$ when these boxes are read from right to
left and then from bottom to top.  We then let $w_{\nu/\lambda} =
s_{i_1} \, s_{i_2} \cdots s_{i_m}$ be the product of the corresponding
simple reflections.  For the skew shape above this gives
$w_{\nu/\lambda} = s_4 \, s_3 \, s_6 \, s_5 \, s_4 \, s_8 \, s_7 \,
s_6 = 1\,2\,5\,7\,3\,9\,4\,6\,8$.  Notice that $w_{\nu/\lambda}$
depends on the numbering of the diagonals as well as the diagram
$\nu/\lambda$.  A theorem of Billey, Jockusch, and Stanley
\cite{billey.jockusch.ea:some} says that a permutation $w$ can be
obtained from a skew diagram in this way, if and only if it is {\em
  321-avoiding}, i.e.\ there are no integers $i < j < k$ for which
$w(i) > w(j) > w(k)$.

Permutations obtained from different numberings of the diagonals in
the same diagram $\nu/\lambda$ differ only by a shift.  In other
words, if $w_{\nu/\lambda}$ is the permutation corresponding to the
numbering which puts the bottom-left box in diagonal number one, then
any other numbering will give a permutation of the form $1^m \times
w_{\nu/\lambda}$.  We can therefore define the stable Grothendieck
polynomial for any skew diagram by $G_{\nu/\lambda}(x;y) =
G_{w_{\nu/\lambda}}(x;y)$.

If the skew diagram is a partition $\lambda$, and if $p$ is the number
of the diagonal containing the box in its upper-left corner, then
$w_\lambda$ is called the {\em Grassmannian permutation for $\lambda$
  with descent in position $p$}.  This permutation is given by
$w_\lambda(i) = i + \lambda_{p+1-i}$ for $1 \leq i \leq p$ and
$w_\lambda(i) < w_\lambda(i+1)$ for $i \neq p$.  Notice that $p$ must
be greater than or equal to the length of $\lambda$, which is the
largest number $\ell = \ell(\lambda)$ for which $\lambda_\ell$ is
non-zero.

It follows from the definitions that the term of lowest total degree
in a Grothen\-dieck polynomial $\Groth_w(x;y)$ is the Schubert
polynomial $\Schub_w(x;-y)$ for the same permutation
\cite{lascoux.schutzenberger:structure}.  This implies that the lowest
term of $G_\lambda(x)$ is the Schur function $s_\lambda(x)$
\cite{macdonald:symmetric*2, macdonald:notes}.  In particular the
polynomials $G_\lambda$ for all partitions $\lambda$ are linearly
independent.  We define $\Gamma$ to be the linear span of all stable
Grothendieck polynomials for Grassmannian permutations:
\[ \Gamma = \bigoplus_\lambda \Z \cdot G_\lambda \subset 
   \Z \llbracket x_1, x_2, \dots, y_1, y_2, \dots \rrbracket \,. 
\]
This group is the main object of study in this paper.  For example,
\refcor{cor_ring} and \refthm{thm_gw} below will show that $\Gamma$ is
closed under multiplication and that it contains all stable
Grothendieck polynomials.

%%% Local Variables: 
%%% mode: latex
%%% TeX-master: "gamma"
%%% End: 

\section{Set-valued tableaux}
\label{sec_settab}

In this section we will introduce set-valued tableaux and use them to
give a formula for stable Grothendieck polynomials indexed by
321-avoiding permutations.

If $a$ and $b$ are two non-empty subsets of the positive integers
$\N$, we will write $a < b$ if $\max(a) < \min(b)$, and $a \leq b$ if
$\max(a) \leq \min(b)$.  We define a {\em set-valued tableau\/} to be
a labeling of the boxes in a Young diagram or skew diagram with finite
non-empty subsets of $\N$, such that the rows are weakly increasing
from left to right and the columns strictly increasing from top to
bottom.  When we speak about a tableau we shall always mean a
set-valued tableau unless explicitly stated otherwise.  The shape
$\sh(T)$ of a tableau $T$ is the partition or skew diagram it is a
labeling of.  For example
\[ \tableau{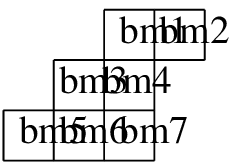} \]
is a tableau of shape $(4,3,3)/(2,1)$ containing the sets $\{2\}$,
$\{3,5\}$, $\{1,2\}$, $\{7\}$, $\{2,3,4\}$, $\{1\}$, and $\{2,3\}$
when the boxes are read bottom to top and then left to right.  A {\em
  semistandard tableau\/} is a set-valued tableau in which each box
contains a single integer.

Given a tableau $T$, let $x^T$ be the monomial in which the exponent
of $x_i$ is the number of boxes in $T$ which contain the integer $i$.
If $T$ is the tableau displayed above we get $x^T = x_1^2\, x_2^4\,
x_3^3\, x_4\, x_5 \, x_7$.  We let $|T|$ denote the total degree of
this monomial, i.e.\ the sum of the cardinalities of the sets in the
boxes of $T$.

\begin{thm}
\label{thm_321}
The single stable Grothendieck polynomial $G_{\nu/\lambda}(x)$ is
given by the formula
\[ G_{\nu/\lambda}(x) = \sum_T (-1)^{|T|-|\nu/\lambda|} \, x^T \]
where the sum is over all set-valued tableaux $T$ of shape
$\nu/\lambda$.
\end{thm}

To prove this proposition we need a result of Fomin and Kirillov.  Let
$H_n(0)$ be the degenerate Hecke algebra over the polynomial ring $R =
\Z[x_1, \dots, x_m]$.  This is the free associative $R$-algebra
generated by symbols $u_1, \dots, u_n$, modulo relations
\begin{align*}
u_i\, u_j &= u_j\, u_i \hspace{49pt} \text{if $|i-j| \geq 2$,} \\
u_i\, u_{i+1}\, u_i &= u_{i+1}\, u_i\, u_{i+1} \hspace{19pt} 
\text{, and} \\
u_i^2 &= - u_i \,.
\end{align*}
If $w \in S_{n+1}$ is a permutation with reduced expression $w =
s_{i_1} \cdots s_{i_\ell}$ we set $u_w = u_{i_1} \cdots u_{i_\ell} \in
H_n(0)$.  This is independent of the choice of a reduced expression.
Furthermore these elements $u_w$ for $w \in S_{n+1}$ form a basis for
$H_n(0)$.  Now set $A(x) = (1 + x\,u_n) \cdots (1 + x\,u_2) \cdot (1 +
x\,u_1)$ and $B(x) = (1 + x\,u_1) \cdot (1 + x\,u_2) \cdots (1 +
x\,u_n)$.  Then a special case of the theory developed by Fomin and
Kirillov \cite{fomin.kirillov:yang-baxter,
  fomin.kirillov:grothendieck} is the following.

\begin{thm}
\label{thm_fk}
The coefficient of $u_w$ in $A(x_m) \cdots A(x_2) \cdot A(x_1)$ is the
stable Groth\-en\-dieck polynomial $G_w(x_1, \dots, x_m)$.  The
coefficient of $u_w$ in $B(x_m) \cdots B(x_2) \cdot B(x_1)$ is $G_w(0;
x_1, \dots, x_m)$.
\end{thm}

Define an {\em inner corner\/} of a partition $\lambda$ to be any box
of $\lambda$ such that the two boxes under and to the right of it are
not in $\lambda$.  By an {\em outer corner\/} we will mean a box
outside $\lambda$ such that the two boxes above and to the left of it
are in $\lambda$.

\begin{proof}[Proof of \refthm{thm_321}]
  It is enough to show that $G_{\nu/\lambda}(x_1,\dots,x_m)$ is the
  sum of the signed monomials $(-1)^{|T|-|\nu/\lambda|} \, x^T$ for
  tableaux $T$ with no integers larger than $m$.
  
  Now number the diagonals of $\nu/\lambda$ from south-west to
  north-east as described in the previous section and let $V =
  \bigoplus_{\mu \supset \lambda} \Z \cdot [\mu]$ be the free Abelian
  group with one basis element for each partition containing
  $\lambda$.  Imitating the methods of Fomin and Greene
  \cite{fomin.greene:noncommutative}, we define a linear action of
  $H_n(0)$ on $V$ as follows.  If a partition $\mu$ has an outer
  corner in the $i$'th diagonal, then we set $u_i \cdot [\mu] =
  [\tilde \mu]$ where $\tilde \mu$ is obtained by adding a box to
  $\mu$ in this corner.  If $\mu$ has an inner corner in the $i$'th
  diagonal, and if the box in this corner is not contained in
  $\lambda$, then we set $u_i \cdot [\mu] = -[\mu]$.  In all other
  cases we set $u_i \cdot [\mu] = 0$.
  
  We claim that if $w \in S_{n+1}$ is any permutation such that $u_w
  \cdot [\lambda] \neq 0$ then $u_w \cdot [\lambda] = [\mu]$ for some
  partition $\mu \supset \lambda$ and $w = w_{\mu/\lambda}$.  This
  claim is clear if $\ell(w) \leq 1$.  If $\ell(w) \geq 2$, write $w =
  s_i \, w'$ where $\ell(w') = \ell(w) - 1$.  Since $u_{w'} \cdot
  [\lambda] \neq 0$ we can assume by induction that $u_{w'} \cdot
  [\lambda] = [\mu_0]$ for some partition $\mu_0$, and $w' =
  w_{\mu_0/\lambda}$.  It is enough to show that $\mu_0$ has no inner
  corner in the $i$'th diagonal.  If it had, then the box in this
  corner would be outside $\lambda$ since $u_i \cdot [\mu_0] \neq 0$.
  But then we could write $w' = w'_1 \, s_i \, w'_2$ such that
  $\ell(w') = \ell(w'_1) + \ell(w'_2) + 1$ and no reduced expression
  for $w'_1$ contains $s_{i-1}$, $s_i$, or $s_{i+1}$.  This would mean
  that $w = s_i \, w' = w'_1 \, w'_2$, contradicting that $\ell(w) >
  \ell(w')$.
  
  Using \refthm{thm_fk} it follows from the claim that the coefficient
  of the basis element $[\nu]$ in $A(x_m) \cdots A(x_1) \cdot
  [\lambda]$ is $G_{\nu/\lambda}(x_1,\dots,x_m)$.  Finally, it is easy
  to identify the terms of $A(x_m) \cdots A(x_1)$ which take
  $[\lambda]$ to $[\nu]$ with set-valued tableaux on $\nu/\lambda$.
  In fact, this product expands as a sum of terms of the form
\[ (x_{i_1} \, u_{j_1}) \cdot (x_{i_2} \, u_{j_2}) \cdots 
   (x_{i_k} \, u_{j_k}) \,.
\]
From any such term which contributes to the coefficient of $[\nu]$ we
obtain a tableau on $\nu/\lambda$ by joining the integer $i_r$ to the
set in the inner corner in diagonal number $j_r$ of the partition
$u_{j_r} \cdots u_{j_k} \cdot [\lambda]$ for each $1 \leq r \leq k$.
\end{proof}

\begin{remark}
  It is easy to extend the notion of set-valued tableau to obtain
  formulas for double stable Grothendieck polynomials
  $G_{\nu/\lambda}(x;y)$.  The main point is that integers
  corresponding to the $x$-variables should occupy horizontal strips
  in a set-valued tableau while integers corresponding to
  $y$-variables should appear in vertical strips.
\end{remark}

While we have the notation of \refthm{thm_fk} fresh in mind, we shall
also establish the following lemma for use later.  At an early point
in this project we asked Sergey Fomin if the statement of this lemma
could possibly be true, after which he proved it immediately.

\begin{lemma}[Fomin]
\label{lemma_conjugate}
Let $w \in S_{n+1}$ be a permutation and let $w_0 \in S_{n+1}$ be the
longest permutation.  Then $G_{w_0 w w_0}(x;y) = G_w(y;x)$.
\end{lemma}
\begin{proof}
Since $G_w(1-e^{-x}; 1-e^y)$ is super-symmetric, it is enough to prove
that $G_{w_0 w w_0}(x_1,\dots,x_m) = G_w(0; x_1,\dots,x_m)$ for any
number of variables $m$.

Now $H_n(0)$ has an $R$-linear automorphism which sends $u_i$ to
$u_{n+1-i}$ for each $1 \leq i \leq n$.  Since this automorphism takes
$A(x_j)$ to $B(x_j)$ and $u_w$ to $u_{w_0 w w_0}$, the lemma follows
from \refthm{thm_fk}.
\end{proof}

As a special case we obtain $G_{\lambda'/\mu'}(x;y) = G_{w_0 \,
  w_{\lambda/\mu} \, w_0} (x;y) = G_{\lambda/\mu}(y;x)$.

%%% Local Variables: 
%%% mode: latex
%%% TeX-master: "gamma"
%%% End: 

\section{A column bumping algorithm}
\label{sec_jdt}

In this section we will present a column bumping algorithm for
set-valued tab\-leaux and derive an important bijective correspondence
from it.  This correspondence will be the main ingredient in the proof
of the Littlewood-Richardson rule for stable Grothendieck polynomials.

We will use the following notation.  If $a$ and $b$ are disjoint sets
of integers we will let \rtab{-7pt}{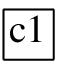} denote a single box containing
the union of $a$ and $b$.  If $T$ is a tableau with $\ell$ columns,
and $C_i$ denotes its $i$'th column for each $1 \leq i \leq \ell$,
then we will write
\[ T = \rtab{-40pt}{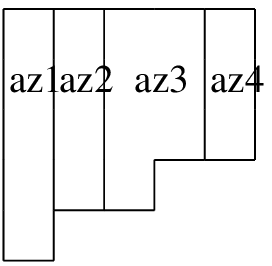} = (C_1, C_2, \dots, C_\ell) \,. \]

We start with the following definition.

\begin{defn}
\label{defn_bump}
Let $x \in \N$, $x_0 \in \N \cup \{\infty\}$, and let $C$ be a
tableau with only one column.  We then define a new tableau $x
\ins{x_0} C$ by the following rules:
\begin{align}
\tag{B1}
 x \xrightarrow[x_0]{} \rtab{-21pt}{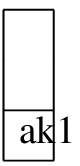} &= \hspace{12pt}
\rtab{-28pt}{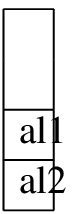} \hspace{1cm}\text{if $a < x$}
\end{align}\begin{align}
%\\
\tag{B2}
 x \xrightarrow[x_0]{} \rtab{-21pt}{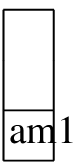} &= \hspace{12pt}
\rtab{-28pt}{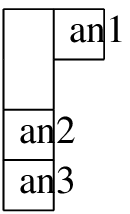} \hspace{.45cm}\text{if $a < x \leq b$}
\\
\tag{B3}
 x \xrightarrow[x_0]{} \rtab{-28pt}{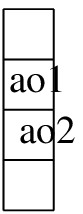} &= \hspace{12pt}
\rtab{-28pt}{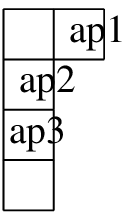} \hspace{.45cm}\text{if $a < x \leq b$}
\\
\tag{B4}
 x \xrightarrow[x_0]{} \rtab{-28pt}{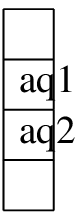} &= \hspace{12pt}
\rtab{-28pt}{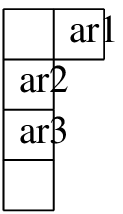} \hspace{.45cm}\text{if $a < x \leq b$ and $x_0 \not \in b$}
\\
\tag{B5}
 x \xrightarrow[x_0]{} \rtab{-21pt}{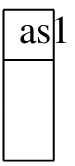} &= \hspace{12pt}
\rtab{-21pt}{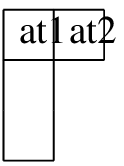} \hspace{.45cm}\text{if $x \leq b$ and $x_0 \not \in b$}
\\
\tag{B6}
 x \xrightarrow[x_0]{} \rtab{-28pt}{B4L} &= \hspace{12pt}
\rtab{-28pt}{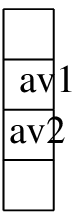} \hspace{1cm}\text{if $a < x \leq b$ and $x_0 \in b$}
\\
\tag{B7}
 x \xrightarrow[x_0]{} \rtab{-21pt}{B5L} &= \hspace{12pt}
\rtab{-21pt}{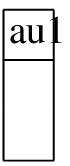} \hspace{1.1cm}\text{if $x \leq b$ and $x_0 \in b$}
\end{align}
\end{defn}
The white areas in these tableaux indicate boxes which are left
unchanged by the operation.  It is easy to see that exactly one of the
cases (B1)--(B7) will apply to define $x \ins{x_0} C$.  In the rules
(B2) to (B5) we say that the set $b$ is ``bumped out''.

If $x \subset \N$ is any non-empty set, we extend this definition as
follows.  Let $x_1 < \dots < x_k$ be the elements of $x$ in increasing
order and let $(C_k, y_k)$ be the tableau $x_k \ins{x_0} C$:
\[ x_k \ins{x_0} C = \rtab{-21pt}{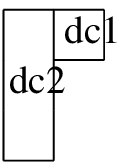} \]
Here $y_k$ is the set in the single box in the second column.  If $x_k
\ins{x_0} C$ has only one column, we let $y_k$ be the empty set.
Continue by setting $(C_i, y_i) = (x_i \ins{x_{i+1}} C_{i+1})$ for
each $i = k-1, \dots, 1$.  We finally define $x \ins{x_0} C$ to be the
tableau $(C_1, y)$ where $y = y_1 \cup \dots \cup y_k$:
\[ x \ins{x_0} C = \rtab{-21pt}{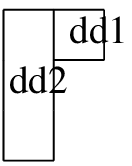} \]
If $x_0 \subset \N$ is a set, we will write $x \ins{x_0} C = x
\ins{\min(x_0)} C$.  We furthermore set $x \cdot C = x \ins{\infty}
C$.  This defines the product of a one-box tableau with a one-column
tableau.  Notice that if $C$ has $\ell$ boxes, then the shape of $x
\cdot C$ is one of $(1^{\ell+1})$, $(2, 1^{\ell-1})$, or $(2,
1^\ell)$, where $1^\ell$ means a sequence of $\ell$ ones.  Notice also
that $x \cdot C = x \ins{x_0} C$ unless one of the rules (B6) or (B7)
are used to define $\max(x) \ins{x_0} C$.  For example:
\[ \rtab{-7pt}{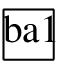} \cdot \rtab{-14pt}{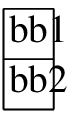} = \rtab{-21pt}{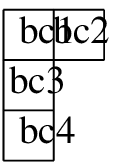} \]

We will continue by defining the product of a non-empty set $x$ with
any tableau $T = (C_1, C_2, \dots, C_\ell)$.  Namely, set $(C_1', y_1)
= x \cdot C_1$ and $(C_i', y_i) = y_{i-1} \cdot C_i$ for $2 \leq i
\leq \ell$.  If a product $y_{i-1} \cdot C_i$ has only one column for
some $i$, we let $C'_i$ be this product and set $C'_j = C_j$ for $j >
i$ and $y_\ell = \emptyset$.  We then define $x \cdot T$ to be the
tableau whose $j$'th column is $C_j'$ for $1 \leq j \leq \ell$.  If
$y_\ell \neq \emptyset$ we furthermore add an $\ell+1$'st column with
one box containing this set:
\[ x ~\cdot~ \rtab{-40pt}{columns} = \hspace{15pt} \rtab{-40pt}{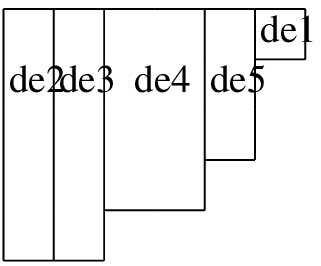} \]
To see that this is in fact a tableau, we need the following lemma.

\begin{lemma}
\label{lemma_setmult}
Let $x \in \N$, $y_0 \in \N \cup \{\infty\}$, and let $C_1$ and
$C_2$ be one-column tableaux which fit together to form a tableau
$(C_1,C_2)$.  Let $(C_1',y) = x \cdot C_1$ and $(C_2',z) = y \ins{y_0}
C_2$.  Then $(C_1', C_2')$ is a tableau.
\end{lemma}
\begin{proof}
  Suppose $y$ was bumped out of box number $i$ in $C_1$, counted from
  the top.  Then since $(C_1,C_2)$ is a tableau, $y$ is less than or
  equal to the set in box number $i$ in $C_2$.  This implies that the
  $j$'th box of $C_2'$ is equal to the $j$'th box of $C_2$ for all $j
  > i$.  For $j \leq i$, the $j$'th box of $C_2'$ can contain elements
  from the $j$'th box of $C_2$ and from $y$.  Since both of these sets
  are larger than or equal to the $j$'th box in $C'_1$, this shows
  that $(C_1',C_2')$ is a tableau.
\end{proof}

To check that the product $x \cdot T$ defined above is a tableau, it
is enough to assume that $T = (C_1, C_2)$ has only two columns.  Let
$x' = \min(x)$ and $x'' = x \smallsetminus \{x'\}$, and set $(C_1'',
y'') = x'' \cdot C_1$ and $(C_2'', z'') = y'' \cdot C_2$.  Then set
$(C_1', y') = x' \ins{x''} C_1''$ and $(C_2', z') = y' \ins{y''}
C_2''$.  Then $x \cdot C_1 = (C_1', y)$ where $y = y'' \cup y'$, and
$y \cdot C_2 = (C_2', z'' \cup z')$.  We must show that $(C_1', C_2')$
is a tableau.

By induction we can assume that $(C_1'', C_2'')$ is a tableau.  If $x'
\ins{x''} C_1''$ is defined by (B6) or (B7) then the maximal
elements in the boxes of $C_1'$ and $C_1''$ are equal, and since $y'$
is empty we have $C_2' = C_2''$. This makes it clear that $(C_1',
C_2')$ is a tableau.  Otherwise we have $x' \ins{x''} C_1'' = x'
\cdot C_1''$ in which case $(C_1', C_2')$ is a tableau by
\reflemma{lemma_setmult}.

Define a {\em rook strip\/} to be a skew shape (between two
partitions) which has at most one box in any row or column.  It then
follows from our earlier observations regarding the shape of a product
of a set with a one-column tableau, that the shape of $x \cdot T$
differs from that of $T$ by a rook strip.

Now let $C$ be a tableau with one column and let $T$ be any tableau.
Suppose the boxes of $C$ contain the sets $x_1, x_2, \dots, x_\ell$,
read from top to bottom.  We then define the product of $C$ and $T$ to
be the tableau $C \cdot T = x_\ell \cdot (x_{\ell-1} \cdot (\cdots x_2
\cdot (x_1 \cdot T) \cdots))$.

\begin{example}
  It is not possible to extend this product to an associative product
  on all set-valued tableau.  In fact, if this was possible we would
  have
\[ \begin{split}
\rtab{-14pt}{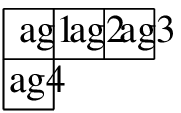}
&= \rtab{-7pt}{2} \cdot \rtab{-14pt}{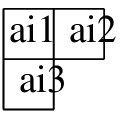}
= \rtab{-7pt}{2} \cdot \rtab{-7pt}{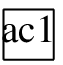} \cdot \rtab{-7pt}{2}
= \rtab{-14pt}{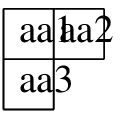} \cdot \rtab{-7pt}{2}
\\
&= \rtab{-14pt}{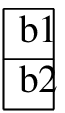} \cdot \rtab{-7pt}{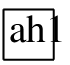} \cdot \rtab{-7pt}{2}
= \rtab{-14pt}{2-1} \cdot \rtab{-14pt}{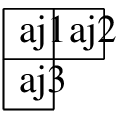}
= \rtab{-14pt}{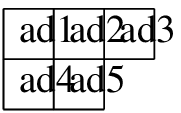}
\end{split} \]
which is of course wrong.
%  It is furthermore possible to show that no
%  associative product exists which reflects multiplication of
%  Grothendieck polynomials.
\end{example}

In the following lemmas we shall study the shape of a product $C \cdot
T$.

\begin{lemma}
\label{lemma_shape1}
Let $x_1 < x_2$ be non-empty sets of integers and let $C$ be a
one-column tableau.  Let $x_1 \cdot C = (C_1, y_1)$ and $x_2 \cdot C_1
= (C_2, y_2)$.  Then $y_1 < y_2$.
\end{lemma}
\begin{proof}
  Notice at first that $\min(x_2) \leq y_2$, which follows directly
  from \refdefn{defn_bump}.  Since all of the integers in $y_2$ come
  from $C_1$, it suffices to show that all integers from $C_1$ which
  are $\geq \min(x_2)$ are also strictly greater than $y_1$.
  
  Let $k$ be the number of the box in $C_1$ (counted from the top)
  which contains $\max(x_1)$.  There are then two possibilities.
  Either $\max(x_1)$ is the largest element in this box of $C_1$, in
  which case all integers in $y_1$ must come from the boxes $1$
  through $k$ of $C$.  Furthermore, the boxes in $C_1$ strictly below
  the $k$'th box contain the same sets as the corresponding boxes of
  $C$, so they are strictly greater than $y_1$.  Since these are also
  the only boxes that can be greater than or equal to $\min(x_2)$, the
  statement follows in this case.
  
  Otherwise $\max(x_1)$ is not the largest element in the $k$'th box
  of $C_1$, which implies that all integers in $y_1$ come from the
  boxes $1$ through $k-1$ in $C$.  Since the integers in $C_1$ which
  are strictly greater than $\max(x_1)$ all come from the $k$'th box
  in $C$ or from boxes below this box, the statement is also true in
  this case.
\end{proof}

\begin{lemma}
\label{lemma_shape2}
Let $T$ be a tableau and let $x_1 < x_2$ be non-empty sets of
integers.  Set $T_1 = x_1 \cdot T$ and $T_2 = x_2 \cdot T_1$, and let
$\theta_1 = \sh(T_1)/\sh(T)$ and $\theta_2 = \sh(T_2)/\sh(T_1)$ be the
rook strips giving the differences between the shapes of these
tableaux.  Then all boxes of $\theta_2$ are strictly south of the
boxes in $\theta_1$.
\end{lemma}
\begin{proof}
  Suppose the southernmost box of $\theta_1$ occurs in column $j$.
  Then let $U$ be the tableau consisting of the leftmost $j-1$ columns
  of $T$ and let $V$ be the rest of $T$.  We will write $T = (U, V)$
  to indicate this.  Similarly we let $(U_1, y_1) = x_1 \cdot U$ and
  $(U_2, y_2) = x_2 \cdot U_1$, i.e.\ $U_1$ and $U_2$ are the leftmost
  $j-1$ columns of these products.  Finally set $V_1 = y_1 \cdot V$
  and $V_2 = y_2 \cdot V_1$.  We then have $T_1 = (U_1, V_1)$ and $T_2
  = (U_2, V_2)$.
  
  Since $V_1$ has one more box in the first column than $V$, this box
  of $V_1$ must contain a subset of $y_1$.  Since $y_2 > y_1$ by
  \reflemma{lemma_shape1}, this means that $V_2 = y_2 \cdot V_1$
  consists of $V_1$ with $y_2$ attached below the first column (or
  $V_2 = V_1$ if $y_2$ is empty).  The lemma follows from this.
\end{proof}

It follows immediately from this lemma that the shape of a product $C
\cdot T$ adds a vertical strip to the shape of $T$, but more detailed
information can be obtained.  As above, let $x_1, x_2, \dots, x_\ell$
be the sets contained in the boxes of $C$ from top to bottom.  Set
$T_0 = T$ and $T_i = x_i \cdot T_{i-1}$ for $1 \leq i \leq \ell$.  Let
$\theta_i = \sh(T_i)/\sh(T_{i-1})$ be the rook strip between the
shapes of $T_i$ and $T_{i-1}$, and let $\theta = \sh(C \cdot
T)/\sh(T)$ be the union of these rook strips.  Then
\reflemma{lemma_shape2} says that the $\theta_i$ split $\theta$ up
into disjoint segments running from north to south.
\[ \theta = \raisebox{-39pt}{\pic{65}{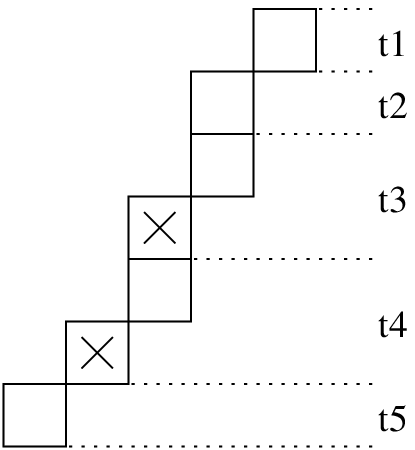}} \]

Define the {\em extra boxes\/} of $\theta$ to be the boxes that are
not the upper-right box of any rook strip $\theta_i$.  These boxes are
marked with a cross in the above picture.  There are exactly $|\theta|
- \ell$ extra boxes in $\theta$, and at most one in each column of
$\theta$.  Furthermore, if a column of $\theta$ has an extra box then
it is the northernmost one.  Define $\col(C, T)$ to be the set of
columns of $\theta$ which contain an extra box.  The rook strips
$\theta_i$ are then uniquely determined by $\theta$ and this set.  In
what follows we shall see that if $x$ is a set and $T$ is a tableau,
then $x$ and $T$ can be recovered from their product $x \cdot T$ if
one knows the shape of $T$.  As a consequence we see that $C$ and $T$
are uniquely determined by their product $C \cdot T$, the shape of
$T$, and the set $\col(C, T)$.

We will state this a bit sharper.  Given a vertical strip $\theta$ and
a non-negative integer $d \geq 0$, let $\CC_d(\theta)$ be the set of
all sets of the non-empty columns of $\theta$ which have cardinality
$d$ and avoids the last column of $\theta$.  If $c(\theta)$ is the
number of non-empty columns of $\theta$ then $\CC_d(\theta)$ has
cardinality $\binom{c(\theta)-1}{d}$.  Let $\TT_\lambda$ be the set of
all set-valued tableaux of shape $\lambda$.  We then have a map
\begin{equation}
\label{eqn_biject} 
  \TT_{(1^\ell)} \times \TT_\lambda  \longrightarrow
  \coprod_\nu \TT_\nu \times \CC_{|\nu/\lambda|-\ell}(\nu/\lambda)
\end{equation}
which takes $(C, T)$ to the pair $(C \cdot T, \col(C, T))$ in the set
$\TT_\nu \times \CC_{|\nu/\lambda|-\ell}(\nu/\lambda)$ where $\nu =
\sh(C \cdot T)$.  The disjoint union is over all partitions $\nu$
containing $\lambda$ such that $\nu/\lambda$ is a vertical strip. 
In the remaining part of this section we will construct an inverse to
this map, thus proving the following.

\begin{thm}
\label{thm_biject}
The map of (\ref{eqn_biject}) is bijective.
\end{thm}

As a first consequence, we obtain a bijective proof of Lenart's Pieri
formula \cite{lenart:combinatorial}.

\begin{cor}[Lenart]
\label{cor_pieri}
For any partition $\lambda$ and $\ell \geq 1$ we have
\[ G_{(1^\ell)} \cdot G_\lambda = \sum_\nu 
   (-1)^{|\nu/\lambda|-\ell}
   \binom{c(\nu/\lambda)-1}{|\nu/\lambda|-\ell} \, G_\nu
\]
where the sum is over all partitions $\nu \supset \lambda$ such that
$\nu/\lambda$ is a vertical strip, and $c(\nu/\lambda)$ is the number
of non-empty columns in this diagram.
\end{cor}

In order to construct the inverse map of (\ref{eqn_biject}) we will
define a reverse column bumping algorithm for set-valued tableaux.  We
start with the following definition.

\begin{defn}
\label{defn_revbump}
Let $T = (C, y)$ be a tableau whose second column has one box
containing a single integer $y \in \N$.  For any $y_0 \in \N \cup
\{0\}$ we define the pair $\RR_{y_0}(C, y)$ by the following rules:
\begin{align}
\tag{R1}
\RR_{y_0}\left( \rtab{-28pt}{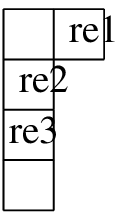} \right) &= 
\left( b ~, ~\rtab{-28pt}{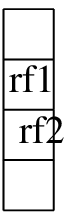} \right) \hspace{15pt}
  \text{if $b \leq y < c$}
\\
\tag{R2}
\RR_{y_0}\left( \rtab{-28pt}{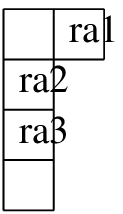} \right) &= 
\left( a ~, ~\rtab{-28pt}{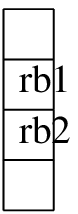} \right) \hspace{15pt}
  \text{if $a \leq y < b$ and $y_0 \not \in a$}
\\
\tag{R3}
\RR_{y_0}\left( \rtab{-21pt}{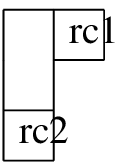} \right) &= 
\left( a ~, ~\rtab{-21pt}{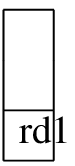} \right) \hspace{15pt}
  \text{if $a \leq y$ and $y_0 \not \in a$}
\\
\tag{R4}
\RR_{y_0}\left( \rtab{-28pt}{R2L} \right) &= 
\left( \emptyset ~, ~\rtab{-28pt}{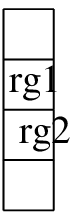} \right) \hspace{15pt}
  \text{if $a \leq y < b$ and $y_0 \in a$}
\\
\tag{R5}
\RR_{y_0}\left( \rtab{-21pt}{R3L} \right) &= 
\left( \emptyset ~, ~\rtab{-21pt}{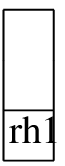} \right) \hspace{15pt}
  \text{if $a \leq y$ and $y_0 \in a$}
\end{align}
\end{defn}

We extend this definition to tableaux $(C, y)$ where $y \subset \N$ is
any non-empty set of integers as follows: Let $y_1 < \dots < y_k$ be
the elements of $y$ in increasing order.  Let $(x_1, C_1) =
\RR_{y_0}(C, y_1)$ and $(x_i, C_i) = \RR_{y_{i-1}}(C_{i-1}, y_i)$ for
$2 \leq i \leq k$.  Then we set $\RR_{y_0}(C, y) = (x, C_k)$ where $x
= x_1 \cup \cdots \cup x_k$.  If $y_0 \subset \N$ is a set, we will
write $\RR_{y_0}(C, y) = \RR_{\max(y_0)}(C, y)$.  Finally we set
$\RR(C, y) = \RR_0(C, y)$.

\begin{lemma}
\label{lemma_inverse}
Let $C$ be a one-column tableau with $\ell$ boxes and let $x \subset
\N$ be a set such that $x \cdot C$ has shape $(2,1^{\ell-1})$.  Then
$\RR(x \cdot C) = (x, C)$.  Similarly, if $T$ is a tableau of shape
$(2, 1^{\ell-1})$ and $\RR(T) = (x, C)$, then $C$ has $\ell$ boxes and
$x \cdot C = T$.
\end{lemma} 
\begin{proof}
  Suppose at first $x \in \N$ is a single integer, and let $x \cdot C
  = (C', y)$.  Then by the definition of $\RR(C', y)$, the minimal
  element of $y$ will bump out $x$ from $C'$, after which the
  remaining elements of $y$ will be added to the same box as the
  minimal element went into.  This recovers the tableau $C$.
  
  If $x$ has more than one element, let $x_1 = \min(x)$ and $x_2 = x
  \smallsetminus \{x_1\}$, and write $x_2 \cdot C = (C_2, y_2)$ and
  $x_1 \ins{x_2} C_2 = (C_1, y_1)$.  Then $x \cdot C = (C_1, y)$ where
  $y = y_1 \cup y_2$.  By induction we can assume $\RR(C_2, y_2) =
  (x_2, C)$.  There are two cases to consider.
  
  Either $y_1 = \emptyset$.  In this case $C_1$ is obtained from $C_2$
  by joining $x_1$ to the box containing $\min(x_2)$.  This means that
  if $\RR(C_2, \min(y_2)) = (x', C')$ then we have $\RR(C_1,
  \min(y_2)) = (\{x_1\} \cup x', C')$.  This proves that $\RR(C_1, y)
  = (x, C)$ in this case.
  
  Otherwise we have $y_1 \neq \emptyset$, in which case $(C_1, y_1) =
  x_1 \cdot C_2$.  We then know that $\RR(C_1, y_1) = (x_1, C_2)$.
  Since $\RR(C_2, y_2) = (x_2, C)$, it is enough to show that
  $\RR(C_2, y_2) = \RR_{y_1}(C_2, y_2)$.  Here it is enough to
  check that none of the rules (R4) or (R5) are used to define
  $\RR_{\max(y_1)}(C_2, \min(y_2))$.  But since $\min(x_2)$ is
  bumped out when $\RR(C_2, \min(y_2))$ is formed, this would imply
  that $\min(x_2)$ and $\max(y_1)$ are in the same box of $C_2$.
  This contradicts the fact that $x_1$ is bumping $y_1$ but not
  $\min(x_2)$ out when multiplied to $C_2$.

  The proof of the second statement is similar and left to the reader.
\end{proof}

To recover the factors of a product $x \cdot C$ which adds two boxes to
the shape of $C$, we continue with the following definition.

\begin{defn}
Let $T = (C, y)$ be a tableau with at least two boxes in the first
column and one box in second column.  Suppose $\RR(C, y)$ has the form
\begin{equation}
\label{eqn_rstar}
\RR(C, y) = \left( x ~,~ \rtab{-28pt}{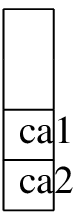} \right) \,. 
\end{equation}
Then we define $\RR^*(C, y)$ as follows:
\begin{align}
\tag{D1}
\RR^*(C,y) &= \left( x \cup b ~,~ \rtab{-21pt}{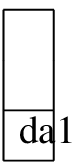} \right)
  \hspace{15pt} \text{if $b \not \subset y$}
\\
\tag{D2}
\RR^*(C,y) &= \left( x ~,~ \rtab{-21pt}{B2L} \right) 
  \hspace{33pt} \text{if $b \subset y$}
\end{align}
\end{defn}

\begin{lemma}
\label{lemma_inverse2}
Let $C$ be a one-column tableau with $\ell$ boxes and let $x \subset
\N$ be a set such that $x \cdot C$ has shape $(2,1^\ell)$.  Then
$\RR^*(x \cdot C) = (x, C)$.  Similarly, if $T$ is a tableau of shape
$(2, 1^\ell)$ and $\RR^*(T) = (x,C)$, then $C$ has $\ell$ boxes and $x
\cdot C = T$.
\end{lemma}
\begin{proof}
  Notice at first that if $x \cdot C$ has shape $(2,1^\ell)$ then
  $\max(x) \cdot C$ must be defined by one of the rules (B1) or (B2).
  Let $x_2 \subset x$ be the largest subset of the form $x \cap
  [k,\infty[$ such that none of the rules (B3)--(B5) are used in the
  definition of $x_2 \cdot C$, and let $x_1 = x \smallsetminus x_2$.
  The lemma is easy to prove if $x_1$ is empty, so we will assume $x_1
  \neq \emptyset$.  Let $x_2 \cdot C = (C_2, y_2)$ and $x_1 \cdot C_2
  = (C_1, y_1)$.  Then $x \cdot C = (C_1, y)$ where $y = y_1 \cup
  y_2$.  Furthermore we have $\RR(C_1, y_1) = (x_1, C_2)$ by
  \reflemma{lemma_inverse}.  There are two cases to consider.

  First suppose $\max(x_2) \cdot C$ is defined by (B1).  Then $C_2$ is
  equal to $C$ with $x_2$ attached in a new box at the bottom, and
  $y_2$ is empty.  This means that $x \cdot C = (C_1, y_1)$, so $\RR(x
  \cdot C) = (x_1, C_2)$.  By rule (D1) we therefore have $\RR^*(x
  \cdot C) = (x_1 \cup x_2, C)$ as claimed.
  
  Otherwise $\max(x_2) \cdot C$ is defined by (B2).  If the bottom box
  of $C$ contains $a \cup b$ with $a < \max(x_2) \leq b$, then $y_2 =
  b$ and $C_2$ is obtained by removing $b$ from the bottom box of $C$
  and attaching $x_2$ in a new box below it.  It follows that
  $\RR_{y_1}(C_2, y_2) = (x_2, C')$ where $C'$ is obtained from $C$ by
  moving the elements of $b$ from the bottom box to a new box below
  it.  We conclude that $\RR(x \cdot C) = (x, C')$, so $\RR^*(x \cdot
  C) = (x, C)$ by (D2).
  
  The proof of the second statement is similar and left to the reader.
\end{proof}

Now let $T$ be a tableau of shape $\nu$ and let $\lambda \varsubsetneq
\nu$ be a proper subpartition such that $\nu/\lambda$ is a rook
strip.  We will then produce a tableau $T'$ of shape $\lambda$ and a set
$x \subset \N$ such that $x \cdot T' = T$.

Write $T = (C_1, C_2, \dots, C_k)$ where $C_i$ is the $i$'th column,
and suppose the upper-right box of $\nu/\lambda$ is in column $j$.
Let $C_j'$ be the result of removing the bottom box from $C_j$, and
let $x_j$ be the set from this removed box.  Now for each $i = j-1,
\dots, 1$, we define
\[ (x_i, C_i') = \begin{cases}
  \RR(C_i, x_{i+1}) & 
    \text{if $\nu/\lambda$ does not have a box in column $i$;} \\
  \RR^*(C_i, x_{i+1}) & 
    \text{if $\nu/\lambda$ has a box in column $i$.}
   \end{cases}
\]
We then set $x = x_1$ and $T' = (C_1', \dots, C_j', C_{j+1}, \dots,
C_k)$.  A argument similar to the proof of \reflemma{lemma_setmult}
shows that $T'$ is a tableau, and by definition the shape of $T'$ is
$\lambda$.  Furthermore, it follows from \reflemma{lemma_inverse} and
\reflemma{lemma_inverse2} that $x \cdot T' = T$.  We let
$\RR_{\nu/\lambda} : \TT_\nu \longrightarrow \TT_1 \times \TT_\lambda$
be the map defined by $\RR_{\nu/\lambda}(T) = (x, T')$.

\begin{proof}[Proof of \refthm{thm_biject}]
  It follows from \reflemma{lemma_inverse} and
  \reflemma{lemma_inverse2} that the maps $\RR_{\nu/\lambda}$ define
  an inverse to the map of (\ref{eqn_biject}) when $\ell = 1$.  If
  $\ell \geq 2$ and $(T, S) \in \TT_\nu \times
  \CC_{|\nu/\lambda|-\ell}(\nu/\lambda)$ is any element, there are
  unique rook strips $\theta_1, \dots, \theta_\ell$ which split the
  vertical strip $\theta = \nu/\lambda$ up into disjoint intervals
  from north to south, such that $S$ contains the columns of the extra
  boxes in $\theta$.  Then set $(x_\ell, T_\ell) =
  \RR_{\theta_\ell}(T)$ and $(x_i, T_i) = \RR_{\theta_i}(T_{i+1})$ for
  $i = \ell-1, \dots, 1$, and let $C \in \TT_{(1^\ell)}$ be the column
  whose $i$'th box contains $x_i$.  An argument similar to the proof
  of \reflemma{lemma_shape1} shows that $x_1 < \dots < x_\ell$ which
  implies that $C$ is a tableau.  Finally \reflemma{lemma_inverse} and
  \reflemma{lemma_inverse2} show that the map $(T, S) \mapsto (C,
  T_1)$ gives an inverse to (\ref{eqn_biject}).
\end{proof}

\begin{remark}
  Although we have skipped some details of the proof of
  \refthm{thm_biject}, the arguments given here do suffice to
  establish that the map of (\ref{eqn_biject}) is injective.  In stead
  of writing down proofs of the remaining statements, one can also use
  Lenart's proof of \refcor{cor_pieri} \cite{lenart:combinatorial} to
  deduce that the two sets in (\ref{eqn_biject}) have the same number
  of elements (which is finite if we only consider tableaux containing
  integers between $1$ and $m$ for any $m \geq 1$).
\end{remark}

%%% Local Variables: 
%%% mode: latex
%%% TeX-master: "gamma"
%%% End: 

\section{Stable Grothendieck polynomials in non-commutative variables}
\label{sec_locplac}

In this section we define stable Grothendieck polynomials in
non-commutative variables and show that they span a commutative ring.
As a consequence we obtain an explicit Littlewood-Richardson rule for
multiplying Grothendieck polynomials of Grassmannian permutations.

Following \cite{fomin:schur, fomin.greene:noncommutative}, we define
the {\em local plactic algebra\/} to be the free associative
$\Z$-algebra $\LL$ in variables $u_1, u_2, \dots$, modulo the
relations
\begin{align*}\,
u_i\, u_j &= u_j\, u_i \hspace{36pt} \text{if $|i-j| \geq 2$,} \\
u_i\, u_{i+1}\, u_i &= u_{i+1}\, u_i\, u_i \hspace{16pt} \text{, and} \\
u_{i+1}\, u_i\, u_{i+1} &= u_{i+1}\, u_{i+1}\, u_i \,.
\end{align*}
We shall work in the completion $\Hat \LL$ of $\LL$ which consists of
formal power series in these variables.

Let $T$ be a set-valued tableau.  We define the (column) {\em word\/}
of $T$ to be the sequence $w(T)$ of the integers contained in its
boxes when these are read from bottom to top and then from left to
right.  The integers within a single box are arranged in increasing
order.  The word of the tableau displayed in the start of
\refsec{sec_settab} is $(2, 3, 5, 1, 2, 7, 2, 3, 4, 1, 2, 3)$.

If $T$ has word $w(T) = (i_1, i_2, \dots, i_\ell)$, then we let $u^T$
be the non-commutative monomial $u^T = u_{i_1} \, u_{i_2} \cdots
u_{i_\ell} \in \LL$.  It is not hard to see that one gets the same
monomial if the boxes of $T$ are read from left to right first and
then from bottom to top, but we shall not need this fact.  If
$\nu/\lambda$ is any skew diagram we then define a stable Grothendieck
polynomial in the variables $u_i$ by
\[ G_{\nu/\lambda}(u) = \sum_T (-1)^{|T|-|\nu/\lambda|} \, u^T  
   \in \Hat \LL
\]
where this sum is over all tableaux $T$ of shape $\nu/\lambda$.

If $C$ is a one-column tableau and $T$ is any tableau, one may easily
verify from \refdefn{defn_bump} that $u^C \cdot u^T = u^{C \cdot T}$
in $\LL$.  From \refthm{thm_biject} we therefore obtain a
non-commutative version of Lenart's Pieri rule.

\begin{lemma}
\label{lemma_locpieri}
If $\lambda$ is any partition and $\ell \geq 1$ an integer, then
\[ G_{(1^\ell)}(u) \cdot G_\lambda(u) = 
   \sum_{\nu} (-1)^{|\nu/\lambda|-\ell}
   \binom{c(\nu/\lambda)-1}{|\nu/\lambda|-\ell} \, G_\nu(u)
\]
holds in $\Hat \LL$, the sum is over all partitions $\nu \supset
\lambda$ such that $\nu/\lambda$ is a vertical strip.
\end{lemma}

A first consequence of this is that the elements $G_{(1^\ell)}(u)$ for
$\ell \geq 1$ commute in $\Hat \LL$.  Now let $\GG = \bigoplus_\lambda
\Z \cdot G_\lambda(u)$ be the span of the polynomials $G_\lambda(u)$
for all partitions $\lambda$, and let $\Hat \GG \subset \Hat \LL$ be
its completion, consisting of all infinite linear combinations of the
$G_\lambda(u)$.

\begin{lemma}
  The group $\Hat \GG$ consists of all formal power series in the
  polynomials $G_{(1^\ell)}(u)$ for $\ell \geq 1$.  In particular
  $\Hat \GG$ is a commutative subring of $\Hat \LL$.
\end{lemma}
\begin{proof}
  Let $n \in \N$ be any positive integer.  It is enough to show that
  for any partition $\lambda$ there exists a polynomial $P_\lambda(u)
  \in \Z \left[ G_{(1^\ell)}(u) \right]_{\ell \geq 1}$ such that
  $G_\lambda(u) - P_\lambda(u)$ is a linear combination of
  Grothendieck polynomials $G_\mu(u)$ for partitions $\mu$ of length
  $\ell(\mu) \geq n$.
  
  Define a partial order on partitions by writing $\mu < \lambda$ if
  $\mu_1 < \lambda_1$, or $\mu_1 = \lambda_1$ and $\ell(\mu) >
  \ell(\lambda)$.  Notice that given a partition $\lambda$ there are
  only finitely many partitions $\mu$ such that $\mu < \lambda$ and
  $\ell(\mu) < n$.
  
  We shall prove that the polynomials $P_\lambda(u)$ exist by
  induction on this order.  Since the smallest partitions $\lambda$
  have only one column, the existence of $P_\lambda(u)$ is clear for
  these partitions.  Let $\lambda = (\lambda_1, \dots, \lambda_\ell)$
  be a partition of length $\ell < n$ and assume that $P_\mu(u)$
  exists for all $\mu < \lambda$.  Let $\sigma = (\lambda_1-1, \dots,
  \lambda_\ell-1)$ be the partition obtained by removing the first
  column of $\lambda$.  Then $\sigma < \lambda$.  By
  \reflemma{lemma_locpieri} we furthermore have
\[ G_\lambda(u) = G_{(1^\ell)}(u) \cdot G_\sigma(u) - \sum_{\nu
  \varsupsetneq \lambda} (-1)^{|\nu/\lambda|-\ell}
  \binom{c(\nu/\lambda)-1}{|\nu/\lambda|-\ell} \, G_\nu(u)
\]
where the sum is over all partitions $\nu$ properly containing
$\lambda$ such that $\nu/\lambda$ is a vertical strip.  Notice that
any such partition $\nu$ for which $|\nu/\lambda| \geq \ell$ must
satisfy $\nu < \lambda$.  We can therefore define
\[ P_\lambda(u) = G_{(1^\ell)}(u) \cdot P_\sigma(u) - \sum_{\nu
  \varsupsetneq \lambda} (-1)^{|\nu/\lambda|-\ell}
  \binom{c(\nu/\lambda)-1}{|\nu/\lambda|-\ell} \, P_\nu(u) \,.
\]
This finishes the proof.
\end{proof}

The lemma shows that any product $G_\lambda(u) \cdot G_\mu(u)$ is an
element of $\Hat \GG$, so we may write
\begin{equation}
\label{eqn_structure}
G_\lambda(u) \cdot G_\mu(u) = \sum_\nu c^\nu_{\lambda \mu} \, G_\nu(u)
\end{equation}
where the coefficients $c^\nu_{\lambda \mu}$ are integers.  Notice
that this linear combination could be infinite.

We say that a sequence of positive integers $w = (i_1, i_2, \dots,
i_\ell)$ is has {\em content\/} $(c_1, c_2, \dots, c_r)$ if $w$
consists of $c_1$ 1's, $c_2$ 2's, and so on up to $c_r$ $r$'s.  If the
content of each subsequence $(i_k, \dots, i_\ell)$ of $w$ is a
partition, then $w$ is called a {\em reverse lattice word}.  Now if
$\nu = (\nu_1, \dots, \nu_p)$ is a partition, we define $u^\nu =
u_p^{\nu_p} \cdots u_2^{\nu_2} \, u_1^{\nu_1} \in \LL$.

\begin{lemma}
\label{lemma_revlat}
  A sequence $w = (i_1, i_2, \dots, i_\ell)$ is a reverse lattice word
  with content $\nu$ if and only if $u_{i_1} u_{i_2} \cdots u_{i_\ell} =
  u^\nu$ in $\LL$.
\end{lemma}
\begin{proof}
  If $w$ is a reverse lattice word, then the rectification of $w$ in
  the plactic monoid is the semistandard Young tableau $U(\nu)$ of
  shape $\nu$ in which all boxes in row $i$ contain the integer $i$
  \cite[Lemma 5.1]{fulton:young}.  This implies that the identity
  $u_{i_1} u_{i_2} \cdots u_{i_\ell} = u^\nu$ holds even with the
  weaker relations of the plactic algebra.
  
  On the other hand, if $u_{i_1} u_{i_2} \cdots u_{i_\ell} = u^\nu$
  then one can obtain the sequence $(p^{\nu_p}, \dots, 2^{\nu_2},
  1^{\nu_1})$ from $w$ by replacing subsequences in the following
  ways:
\begin{align*}
(i,j) &\leftrightarrow (j,i) \hspace{35pt} \text{if $|i-j| \geq 2$} \\
(i,i+1,i) &\leftrightarrow (i+1,i,i) \\
(i+1,i,i+1) &\leftrightarrow (i+1,i+1,i)
\end{align*}
Since all of these moves preserve reverse lattice words, $w$ must be a
reverse lattice word with content $\nu$.
\end{proof}

If $\lambda$ and $\mu$ are partitions, we let $\lambda * \mu$ be the
skew diagram obtained by putting $\lambda$ and $\mu$ corner to corner
as shown:
\[ \lambda * \mu ~=~ 
   \raisebox{-32pt}{\includegraphics[scale=0.7]{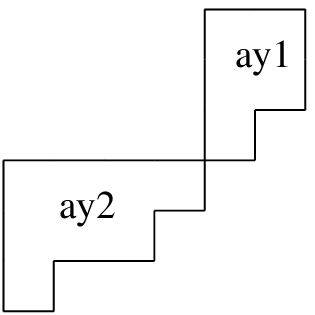}} 
\]

\begin{thm}
\label{thm_lrmult}
The coefficient $c^\nu_{\lambda \mu}$ is equal to
$(-1)^{|\nu|-|\lambda|-|\mu|}$ times the number of set-valued tableaux
$T$ of shape $\lambda * \mu$ such that $w(T)$ is a reverse lattice
word with content $\nu$.
\end{thm}
\begin{proof}
  Start by noticing that the only tableau of shape $\nu$ whose word is
  a reverse lattice word is the tableau $U(\nu)$.  It follows from
  this that the coefficient of $u^\nu$ in the right hand side of
  (\ref{eqn_structure}) is $c^\nu_{\lambda \mu}$.
  
  On the other hand, the left hand side is equal to $G_{\lambda *
    \mu}(u)$.  If $T$ is a tableau on $\lambda * \mu$, then $u^T$ is
  equal to $u^\nu$ exactly when $w(T)$ is a reverse lattice word with
  content $\nu$ by the lemma.  The theorem follows from this.
\end{proof}

Since there are only finitely many tableaux $T$ of a given shape such
that the word of $T$ is a reverse lattice word, \refthm{thm_lrmult}
implies that the linear combination in (\ref{eqn_structure}) is
finite.  In other words $\GG$ is a commutative subring of $\Hat \LL$.

\begin{cor}
\label{cor_ring}
  $\Gamma = \bigoplus \Z \cdot G_\lambda$ is closed under
  multiplication.  The structure constants $c^\nu_{\lambda \mu}$ such
  that $G_\lambda \cdot G_\mu = \sum_\nu c^\nu_{\lambda \mu}\, G_\nu$
  are given by \refthm{thm_lrmult}.
\end{cor}
\begin{proof}
  By replacing each non-commutative variable $u_i$ with $x_i$ in
  \refeqn{eqn_structure}, we obtain the identity $G_\lambda(x) \cdot
  G_\mu(x) = \sum_\nu c^\nu_{\lambda \mu} \, G_\nu(x)$ for single
  stable Grothendieck polynomials.  Since $G_w(1-e^{-x}; 1-e^y)$ is
  super symmetric, the same equation must hold for the double
  polynomials as well.
\end{proof}

\begin{example}
For the shape $(1) * (1)$ we can find the following three set-valued
tableaux whose contents are reverse lattice words:
\[ \tableau{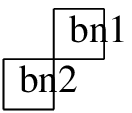} \hspace{20pt} 
   \tableau{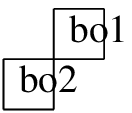} \hspace{20pt} \tableau{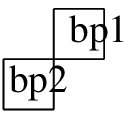} 
\]
It follows that $G_1 \cdot G_1 = G_2 + G_{(1,1)} - G_{(2,1)}$.
\end{example}

Our methods seem insufficient to prove that a stable Grothendieck
polynomial $G_{\nu/\lambda}(u)$ in non-commutative variables is in
$\Hat \GG$, except when the skew shape $\nu/\lambda$ is a product of
partitions like in \refthm{thm_lrmult}.  If this could be established,
then the proof of \refthm{thm_lrmult} would also prove a rule for
writing $G_{\nu/\lambda}$ as a linear combination of the stable
polynomials $G_\mu$ for Grassmannian permutations.  We shall instead
derive such a formula from the statement of \refthm{thm_lrmult} in the
next section.

%%% Local Variables: 
%%% mode: latex
%%% TeX-master: "gamma"
%%% End: 

\section{A coproduct on stable Grothendieck polynomials}
\label{sec_coprod}

Our main task in this section is to show that the ring $\Gamma$ has a
natural coproduct $\Delta : \Gamma \to \Gamma \otimes \Gamma$ which
makes it a bialgebra.  We will show that for certain integers
$d^\nu_{\lambda \mu}$ given by an explicit Littlewood-Richardson rule
similar to that of \refthm{thm_lrmult}, we have
\begin{equation}
\label{eqn_coprod}
  G_\nu(x,z; w,y) = \sum_{\lambda,\mu} d^\nu_{\lambda \mu} \,
   G_\lambda(x;w) \cdot G_\mu(z;y)
\end{equation}
whenever $x$, $y$, $z$, and $w$ are different sets of variables.  The
coproduct can then be defined by $\Delta \, G_\nu = \sum
d^\nu_{\lambda \mu}\, G_\lambda \otimes G_\mu$.

\begin{lemma}
\label{lemma_split1}
  Let $\lambda \subset \nu$ be partitions and let $p, q \geq 1$ be
  integers.  Then
\[ 
   G_{\nu/\lambda}(x_1,\dots,x_{p+q}) =
   \sum_{\tau/\sigma \text{ \rm rook strip}}
   (-1)^{|\tau/\sigma|} \, G_{\tau/\lambda}(x_1,\dots,x_p) \cdot
   G_{\nu/\sigma}(x_{p+1},\dots,x_{p+q})
\]
where the sum is over all partitions $\sigma$ and $\tau$ such that
$\lambda \subset \sigma \subset \tau \subset \nu$ and $\tau/\sigma$
is a rook strip.
\end{lemma}
\begin{proof}
  By \refthm{thm_321} the left hand side is the signed sum of
  monomials $x^T$ for all tableaux $T$ of shape $\nu/\lambda$ such
  that all contained integers are between $1$ and $p+q$.  If $T$ is
  such a tableaux, let $T_1$ be the subtableau obtained by removing
  all integers strictly larger than $p$ (as well as all boxes that
  become empty), and let $T_2$ be obtained by removing integers less
  than or equal to $p$.  Then the shape of $T_1$ is $\tau/\lambda$ for
  some partition $\tau$, while the shape of $T_2$ is of the form
  $\nu/\sigma$, and we have $\lambda \subset \sigma \subset \tau
  \subset \nu$.  Since $\tau/\sigma$ is the shape where $T_1$ and
  $T_2$ overlap, this skew shape must be a rook strip.  Finally, since
  $x^T = x^{T_1} \cdot x^{T_2}$ this gives a bijection between the
  terms of the two sides of the claimed identity.
\end{proof}

Notice that since the polynomials $G_{\nu/\lambda}(1-e^{-x}; 1-e^y)$ are
super-symmetric, this lemma implies that if $x$, $y$, $z$, and $w$ are
different set of variables, then we have
\begin{equation}
  G_{\nu/\lambda}(x,z; w,y) = \sum_{\tau/\sigma \text{ rook strip}}
  (-1)^{|\tau/\sigma|} \, G_{\tau/\lambda}(x;w) \cdot 
  G_{\nu/\sigma}(z;y) \,.
\end{equation}
This can be deduced by writing each polynomial $G_w(1-e^{-x};1-e^y)$
as a linear combination of double Schur functions
\cite{macdonald:symmetric*2}.

\begin{lemma}
\label{lemma_cutshape}
  Let $\theta$ be a skew shape which is broken up into two smaller
  skew shapes $\theta_1$ and $\theta_2$ by a vertical line as shown.
  \[ \theta = \raisebox{-31.5pt}{\pic{53}{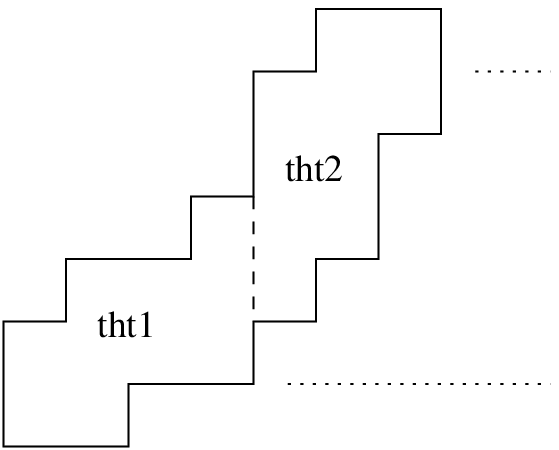}} 
     \hspace{-30cm} \left\lbrace \hspace{30cm}
     \begin{matrix} \vspace{35pt} \\ \end{matrix} \right\rbrace p
  \]
  Let $p$ be the number of boxes between the top edge of the leftmost
  column of $\theta_2$ and the bottom edge of the rightmost column of
  $\theta_1$.  Then we have $G_\theta(x_1,\dots,x_p) =
  G_{\theta_1}(x_1,\dots,x_p) \cdot G_{\theta_2}(x_1,\dots,x_p)$.
\end{lemma}
\begin{proof}
  Number the rows of $\theta$ such that the top box in the leftmost
  column of $\theta_2$ is in row number one, and so that the numbers
  increase from top to bottom.  Then suppose $T_1$ and $T_2$ are
  tableaux of shapes $\theta_1$ and $\theta_2$ for which all contained
  integers are less than or equal to $p$.  Then all integers in row
  $i$ of $T_2$ will be greater than or equal to $i$ because they have
  at least $i-1$ boxes above them.  Similarly the integers in row $i$
  of $T_1$ will be smaller than or equal to $i$ because they have
  $p-i$ boxes below them.  Therefore $T_1$ and $T_2$ fit together to
  form a tableau $T$ of shape $\theta$.  This shows that the terms of
  $G_\theta$ and $G_{\theta_1} \cdot G_{\theta_2}$ are in bijective
  correspondence.
\end{proof}

\begin{prop}
  \label{prop_cutpart}
  Let $\nu$ be a Young diagram which is broken up into two smaller
  Young diagrams $\lambda$ and $\mu$ by a vertical line after
  column $q$.
\[ \nu = \raisebox{-25pt}{\pic{53}{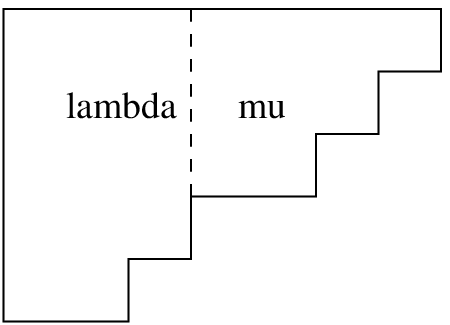}} \]
Then if $p$ is the length of the last column of $\lambda$ we have
\[ G_\nu(x_1,\dots,x_p; y_1,\dots,y_q) =
   G_\lambda(x_1,\dots,x_p; y_1,\dots,y_q) \cdot 
   G_\mu(x_1,\dots,x_p) \,.
\]
\end{prop}
\begin{proof}
  In this proof we will write $x$ for the variables $x_1,\dots,x_p$
  and $y$ for $y_1,\dots,y_q$.  Then by \reflemma{lemma_split1} we
  have
\[ G_\nu(x; y) = \sum G_\tau(0; y) \cdot G_{\nu/\sigma}(x) \]
where the sum is over all partitions $\sigma \subset \tau \subset \nu$
such that $\tau/\sigma$ is a rook strip.  Notice that when $\tau$ has
more than $q$ columns, then $G_\tau(0;y) = G_{\tau'}(y) = 0$ by
\reflemma{lemma_conjugate} and \refthm{thm_321}.  Therefore we only
need to include terms for which $\tau \subset \lambda$ in the sum.
For such terms \reflemma{lemma_cutshape} implies that
$G_{\nu/\sigma}(x) = G_{\lambda/\sigma}(x) \cdot G_\mu(x)$.  The lemma
follows from this by applying \reflemma{lemma_split1} to $G_\lambda(x;
y)$.
\end{proof}

\begin{cor}
\label{cor_zero}
If $\nu_{p+1} \geq q+1$ then $G_\nu(x_1,\dots,x_p; y_1,\dots,y_q) = 0$.
\end{cor}
\begin{proof}
Let $\lambda$ be the first $q$ columns of $\nu$ and let $\mu$ be rest
like in \refprop{prop_cutpart}.  Then since $\ell(\mu) > p$ we get
$G_\mu(x_1,\dots,x_p) = 0$.  The statement therefore follows from the
proposition.
\end{proof}

\begin{cor}[Factorization formula]
\label{cor_factor}
Let $R = (q)^p$ be a rectangle with $p$ rows and $q$ columns, and let
$\sigma$ and $\tau$ be partitions such that $\ell(\sigma) \leq p$ and
$\tau_1 \leq q$.  Let $(R+\sigma,\tau)$ denote the partition
$(q+\sigma_1, \dots, q+\sigma_p, \tau_1, \tau_2, \dots)$ obtained by
attaching $\sigma$ and $\tau$ to the right and bottom sides of $R$.
Then
\begin{multline*}
G_{R+\sigma,\tau}(x_1,\dots,x_p; y_1,\dots,y_q) = \\
   G_\tau(0; y_1,\dots,y_q) \cdot G_R(x_1,\dots,x_p; y_1,\dots,y_q)
   \cdot G_\sigma(x_1,\dots,x_p) \,.
\end{multline*}
\end{cor}
\begin{proof}
Let $x$ denote $x_1,\dots,x_p$ and let $y$ denote $y_1,\dots,y_q$.
Then \refprop{prop_cutpart} implies that $G_{R+\sigma,\tau}(x;y) =
G_{R,\tau}(x;y) \cdot G_\sigma(x)$.  Now using
\reflemma{lemma_conjugate} we obtain $G_{R,\tau}(x;y) =
G_{R'+\tau'}(y;x) = G_{R'}(y;x) \cdot G_{\tau'}(y) =
G_\tau(0;y) \cdot G_R(x;y)$.
\end{proof}

We are now ready to define the structure constants for the coproduct
in $\Gamma$.  We will say that a sequence of integers $w = (i_1,
\dots, i_\ell)$ is a {\em partial reverse lattice word\/} with respect
to an integer interval $[a,b]$, if for all $1 \leq k \leq \ell$ and $p
\in [a,b-1]$, the subsequence $(i_k, \dots, i_\ell)$ has more
occurrences of $p$ than of $p+1$.

Given three partitions $\lambda = (\lambda_1,\dots,\lambda_p)$, $\mu =
(\mu_1,\dots,\mu_q)$, and $\nu$, we then define $d^\nu_{\lambda \mu}$
to be $(-1)^{|\lambda|+|\mu|-|\nu|}$ times the number of set-valued
tableaux $T$ of shape $\nu$ such that $w(T)$ is a partial reverse
lattice word with respect to both of the intervals $[1,p]$ and
$[p+1,p+q]$, and with content $(\lambda,\mu) =
(\lambda_1,\dots,\lambda_p,\mu_1,\dots,\mu_q)$.  Notice that if $R$ is
a rectangle which is taller than $\lambda$ and wider than $\mu$ then
$d^\nu_{\lambda \mu} = c^{R+\lambda,\mu}_{\nu,R}$ by
\refthm{thm_lrmult}.

\begin{thm}
\label{thm_lrsplit}
For any partition $\nu$ we have
\[ G_\nu(x; y) = \sum_{\lambda, \mu} d^\nu_{\lambda \mu} \, 
   G_\lambda(x) \cdot G_\mu(0; y) \,.
\]
\end{thm}
\begin{proof}
  It is enough to show this for finitely many variables
  $x_1,\dots,x_p$ and $y_1,\dots,y_q$, as long as $p$ and $q$ can be
  arbitrarily large.  Let $R = (q)^p$ be a rectangle with $p$ rows and
  $q$ columns.  If $\rho$ is a partition such that $G_\rho$ occurs
  with non-zero coefficient in $G_R \cdot G_\nu$ then first of all $R
  \subset \rho$.  Furthermore, \refcor{cor_zero} shows that
  $G_\rho(x;y)$ is non-zero only if $\rho$ has the form $\rho =
  (R+\lambda,\mu)$ for partitions $\lambda$ and $\mu$.  By these
  observations we get
\[ G_R(x;y) \cdot G_\nu(x;y) = 
   \sum_{\lambda,\mu} d^\nu_{\lambda \mu} \, G_{R+\lambda,\mu}(x;y) =
   \sum_{\lambda,\mu} d^\nu_{\lambda \mu} \, G_R(x;y) \cdot
   G_\lambda(x) \cdot G_\mu(0;y) \,.
\]
Since \refthm{thm_321} shows that $G_R(x;y) \neq 0$, this proves the
theorem.
\end{proof}

Again using the fact that $G_\lambda(1-e^{-x}; 1-e^y)$ is
super-symmetric, this theorem implies that $G_\nu(x,z; w,y) = \sum
d^\nu_{\lambda \mu} \, G_\lambda(x;w) \cdot G_\mu(z;y)$ whenever $x$,
$y$, $z$, and $w$ are different sets of variables.

\begin{cor}
\label{cor_Gamma}
The ring $\Gamma = \bigoplus \Z \cdot G_\lambda$ is a commutative and
cocommutative bialgebra with product $\Gamma \otimes \Gamma \to
\Gamma$ given by
\[ G_\lambda \cdot G_\mu = \sum_\nu c^\nu_{\lambda \mu} \, G_\nu \]
and coproduct $\Delta : \Gamma \to \Gamma \otimes \Gamma$ given by
\[ \Delta \, G_\nu = \sum_{\lambda,\mu} d^\nu_{\lambda \mu} \,
   G_\lambda \otimes G_\mu \,.
\]
The linear map $\Gamma \to \Z$ sending $G_{(0)} = 1$ to one and
$G_\lambda$ to zero for $\lambda \neq \emptyset$ is a counit.
Furthermore, conjugation of partitions defines an involution of this
bialgebra
\[ G_\lambda(x;y) \mapsto G_{\lambda'}(x;y) = G_\lambda(y;x) \,. \]
\end{cor}
\begin{proof}
These statements are clear from Theorems \ref{thm_lrmult} and
\ref{thm_lrsplit}, and \reflemma{lemma_conjugate}.
\end{proof}

\begin{example}
\label{exm_splitone}
  Using the definition of the coefficients $d^\nu_{\lambda \mu}$ we
  may compute
\[ \Delta \, G_1 = G_1 \otimes 1 + 1 \otimes G_1 - G_1 \otimes G_1 \,. 
\]
\end{example}

We will finish this section by showing that $\Gamma$ indeed contains
all stable Grothen\-dieck polynomials.  We start by generalizing the
Littlewood-Richardson rule of \refthm{thm_lrmult} to a rule for the
stable Grothendieck polynomial of any 321-avoiding permutation.

Given partitions $\lambda \subset \nu$ we define $G_{\nu \sslash
  \lambda} = \sum_\mu d^\nu_{\lambda \mu} \, G_\mu$.  It then follows
from \refcor{cor_Gamma} that 
\begin{equation}
\label{eqn_defsslash}
  \Delta \, G_\nu = \sum_{\lambda \subset \nu} 
  G_\lambda \otimes G_{\nu \sslash \lambda} \,,
\end{equation}
and $G_{\nu \sslash \lambda}$ is furthermore uniquely defined by this
identity.  Therefore we deduce from \reflemma{lemma_split1} that
\begin{equation}
\label{eqn_sslash}
  G_{\nu \sslash \lambda} = \sum_{\lambda/\sigma \text{ rook strip}}
  (-1)^{|\lambda/\sigma|} \, G_{\nu/\sigma}
\end{equation}
where the sum is over all partitions $\sigma \subset \lambda$ such
that $\lambda/\sigma$ is a rook strip.

Now given a skew shape $\theta = \nu/\lambda$ and a partition $\mu$,
let $\alpha_{\theta,\mu} = \alpha_{\nu/\lambda,\mu}$ be
$(-1)^{|\mu|-|\theta|}$ times the number of set-valued tableaux $T$ of
shape $\theta$ such that $w(T)$ is a reverse lattice word with content
$\mu$.

\begin{thm}
\label{thm_lr321}
For any skew shape $\theta = \nu/\lambda$ we have
\[ G_{\nu/\lambda} = \sum_{\mu} \alpha_{\theta,\mu} \, G_\mu \,. \]
\end{thm}
\begin{proof}
  We will start by comparing the coefficients $d^\nu_{\lambda \mu}$
  and $\alpha_{\theta,\mu}$.  Suppose $T$ is a tableau of shape $\nu$
  such that $w(T)$ is a partial reverse lattice word with respect to
  the intervals $[1, \ell(\lambda)]$ and $[\ell(\lambda)+1,
  \ell(\lambda) + \ell(\mu)]$ and with content $(\lambda,\mu)$.  Then
  all integers in $T$ which come from the interval $[1,\ell(\lambda)]$
  must be located in the upper-left corner in $T$ of shape $\lambda$.
  Furthermore, any such integer $i$ can occur only in row $i$.  Now
  let the skew shape $\nu/\sigma$ be the region in which the integers
  larger than $\ell(\lambda)$ are located in $T$.  Since this region
  can only overlap $\lambda$ in a rook strip, $\lambda/\sigma$ must be
  a rook strip.  If you remove all integers smaller than or equal to
  $\ell(\lambda)$ from $T$ and subtract $\ell(\lambda)$ from the rest,
  then the result is a tableau of shape $\nu/\sigma$ whose word is a
  reverse lattice word with content $\mu$.  Since $T$ is uniquely
  determined by $\lambda$ and this skew tableau, we obtain
\[ d^\nu_{\lambda \mu} = \sum_{\lambda/\sigma \text{ rook strip}}
   (-1)^{|\lambda/\sigma|} \, \alpha_{\nu/\sigma,\mu}
\]
where this sum is over all partitions $\sigma \subset \lambda$ such
that $\lambda/\sigma$ is a rook strip.

To finish the proof we set $\Tilde G_{\nu/\sigma} = \sum_\mu
\alpha_{\nu/\sigma, \mu} \, G_\mu$.  Then we obtain
\[ G_{\nu \sslash \lambda} =
   \sum_\mu \sum_{\lambda/\sigma \text{ rook strip}} 
     (-1)^{|\lambda/\sigma|} \, \alpha_{\nu/\sigma, \mu} \, G_\mu = 
   \sum_{\lambda/\sigma \text{ rook strip}} (-1)^{|\lambda/\sigma|} \,
     \Tilde G_{\nu/\sigma} \,.
\]
By comparing with (\ref{eqn_sslash}) and noting that the transition
matrix between the $G_{\nu/\lambda}$ and the $G_{\nu \sslash \lambda}$
for fixed $\nu$ is invertible, we conclude that $G_{\nu/\lambda} =
\Tilde G_{\nu/\lambda}$.
\end{proof}

Let us remark that all of the theorems \ref{thm_321},
\ref{thm_lrmult}, \ref{thm_lrsplit}, and \ref{thm_lr321} can be
summarized in the following statement.  We leave the details to the
reader.

\begin{cor}
The $n-1$ fold coproduct applied to $G_{\nu/\lambda}$ is given by
\[ \Delta^{n-1} \, G_{\nu/\lambda} = \sum_{\mu(1),\dots,\mu(n)}
   d^{\nu/\lambda}_{\mu(1),\dots,\mu(n)} \,
   G_{\mu(1)} \otimes \cdots \otimes G_{\mu(n)}
\]
where $d^{\nu/\lambda}_{\mu(1),\dots,\mu(n)}$ is
$(-1)^{|\nu/\lambda|+\sum |\mu(i)|}$ times the number of set-valued
tableaux $T$ of shape $\nu/\lambda$ such that $w(T)$ is a partial
reverse lattice word with respect to each of the intervals $ \left[ 1
  + \sum_{i=1}^{k-1} \ell(\mu(i)) ~,~ \sum_{i=1}^k \ell(\mu(i))
\right]$ for $1 \leq k \leq n$, and has content $(\mu(1), \dots,
\mu(n))$.
\end{cor}

Finally, we will give an independent argument showing that the stable
Grothen\-dieck polynomial $G_w = G_w(x;y)$ for any permutation $w$ is
contained in $\Gamma$.  Recall that the single Schur function
$s_\lambda(x)$ for a partition $\lambda$ is defined by
\[ s_\lambda(x) = \sum_T x^T \]
where the sum is over all semistandard tableaux $T$ of shape $\lambda$
\cite{macdonald:symmetric*2, fulton:young}.  This is the term of
lowest degree in the single stable polynomial $G_\lambda(x)$.

If $\mu$ is a partition containing $\lambda$, let $g_{\lambda \mu}$ be
the number of row and column strict tableaux of shape $\mu/\lambda$
such that all entries in the $i$'th row are between $1$ and $i-1$.
The boxes of these tableaux should contain single integers.  Then
Lenart has proved the following result \cite{lenart:combinatorial}.

\begin{thm}[Lenart]
\label{thm_gtos}
Let $w_\lambda$ be the Grassmannian permutation for $\lambda$ with
descent in position $p \geq \ell(\lambda)$.  Then the single
Grothendieck polynomial $\Groth_{w_\lambda}(x)$ is given by
\[ \Groth_{w_\lambda}(x) = \sum_{\mu \supset \lambda} 
   (-1)^{|\mu/\lambda|} g_{\lambda \mu} \, s_\mu(x_1,\dots,x_p)
\]
where the sum is over all partitions $\mu$ containing $\lambda$.
\end{thm}

As a consequence we obtain the formula
\begin{equation}
\label{eqn_gtos}
  G_\lambda(x) = \sum_{\mu \supset \lambda} (-1)^{|\mu/\lambda|}
  g_{\lambda \mu} \, s_\mu(x) \,.
\end{equation}
This formula can also be derived from \refthm{thm_321}.  We
will briefly sketch the argument.  It is sufficient to construct a
bijection between set-valued tableaux $T$ of shape $\lambda$ and pairs
$(U,T)$ where $U$ is a semistandard tableau of some shape $\mu$
containing $\lambda$ and $S$ is one of the row and column strict
tableaux contributing to $g_{\lambda \mu}$.  We shall do this by
induction on $\ell = \ell(\lambda)$.

Given a set-valued tableau $T$ of shape $\lambda$, let $R$ be the top
row of $T$ and let $T'$ be the rest.  By induction we can assume that
$T'$ corresponds to a pair $(U', S')$.  Now let $\Tilde R$ be the
unique semistandard tableau of shape $(\lambda_1, 1^m)$ where $m = |R|
- \lambda_1$, such that $x^{\Tilde R} = x^R$ and each box in the top
row of $\Tilde R$ contains the smallest integer in the corresponding
box of $R$.  For example:
\[ R = \rtab{-7pt}{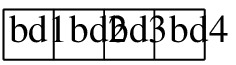} ~~\text{ gives }~~ 
   \Tilde R = \rtab{-28pt}{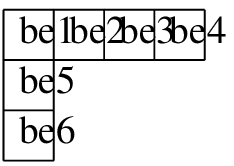}
\]
Now let $U = U' \cdot \Tilde R$ be the product of $U'$ and $\Tilde R$
in the sense of \cite[\S 1]{fulton:young}.  Then $U$ has shape $\mu =
(\lambda_1, \sigma)$ where the partition $\sigma$ is obtained by
adding a vertical strip to the shape of $U'$.  Finally, let $S$ be the
skew tableau of shape $\mu/\lambda$ obtained by copying the $i$'th row
of $S'$ to the $i+1$'st row of $S$; if $S$ needs an additional box in
this row, we put $i$ in this box.  One may then check that the map $T
\mapsto (U,S)$ gives the desired bijection.  For example:
\[ T = \rtab{-14pt}{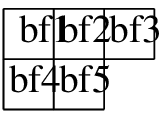} ~~\text{ gives }~~ (U,S) = 
   \left( \rtab{-28pt}{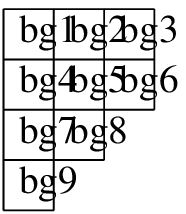} ~,~ \rtab{-28pt}{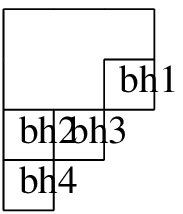} \right) 
\]

Equation (\ref{eqn_gtos}) shows that a stable polynomial
$G_\lambda(x)$ is an infinite linear combination of Schur functions
$s_\mu(x)$ for partitions $\mu \supset \lambda$ such that $\mu$ and
$\lambda$ have the same number of columns.  Since the coefficient of
$s_\lambda(x)$ is $g_{\lambda \lambda} = 1$, we can also go in the
opposite direction and write each Schur function $s_\mu(x)$ as an
infinite linear combination of single stable Grothendieck polynomials
for partitions with the same number of columns as $\mu$.  In fact
Lenart also gives an explicit formula for doing this
\cite{lenart:combinatorial}.

We also need the following result of Fomin and Greene
\cite{fomin.greene:noncommutative}.  If $w \in S_n$ is a permutation
and $\lambda$ is a partition, let $g_{w \lambda}$ be the number of
semistandard tableaux $T$ of shape $\lambda'$ such that if the column
word of $T$ is $w(T) = (i_1, \dots, i_r)$ then $u_{i_1} \cdots u_{i_r}
= \pm u_w$ in the degenerate Hecke algebra $H_n(0)$.  Notice that this
implies that $g_{w \lambda}$ is zero if $\lambda$ has $n$ or more
columns.

\begin{thm}[Fomin and Greene]
\label{thm_fg}
The single stable Grothendieck polynomial $G_w(x)$ is given by
\[ G_w(x) = \sum_\lambda (-1)^{|\lambda|-\ell(w)} g_{w \lambda} \,
   s_\lambda(x) \,. 
\]
\end{thm}

We can finally prove:

\begin{thm}
\label{thm_gw}
Let $w \in S_n$ be a permutation.  Then the double stable Grothendieck
polynomial $G_w(x;y)$ can be written as a finite linear combination
\begin{equation}
\label{eqn_alpha}
  G_w(x;y) = \sum_\lambda \alpha_{w \lambda} \, G_\lambda(x;y) 
\end{equation}
of the polynomials $G_\lambda(x;y)$ for $|\lambda| \geq \ell(w)$.  In
particular $G_w$ is an element of $\Gamma$.
\end{thm}
\begin{proof}
  It follows from (\ref{eqn_gtos}) and \refthm{thm_fg} that the single
  polynomial $G_w(x)$ is a possibly infinite linear combination of
  polynomials $G_\lambda(x)$ for partitions $\lambda$ with at most
  $n-1$ columns:
  \[ G_w(x) = \sum_\lambda \alpha_{w \lambda} G_\lambda(x)\,. \]
  Using that $G_w(1-e^{-x}; 1-e^y)$ is super symmetric, this implies
  that (\ref{eqn_alpha}) is true with the same coefficients.
  
  Let $w_0 = n \cdots 2 \, 1$ be the longest permutation of $S_n$.
  Using \reflemma{lemma_conjugate} we then deduce that $\alpha_{w
    \lambda} = \alpha_{w_0 w w_0, \lambda'}$ is zero unless
  $\ell(\lambda) = \lambda'_1 < n$.  We conclude that $G_w(x;y)$ is a
  finite linear combination of the polynomials $G_\lambda(x;y)$.
\end{proof}

Since the term of lowest degree in $G_w(x)$ is the stable Schubert
polynomial or Stanley symmetric function $F_w(x)$, it follows that 
when $|\lambda| = \ell(w)$ the coefficient $\alpha_{w \lambda}$ is the
one defined by Stanley \cite{stanley:on*1}.  The coefficients
$\alpha_{w \lambda}$ also generalize the structure constants
$c^\nu_{\lambda \mu}$ and $d^\nu_{\lambda \mu}$ of $\Gamma$.  On the
other hand the coefficients $\alpha_{w \lambda}$ are special cases
of quiver coefficients $c_r(\mu)$ which will be introduced in
\cite{buch:structure}.  We believe that these quiver coefficients have
alternating signs.  A consequence of this is the following.

\begin{conj}
\label{conj_alpha}
  The coefficients $\alpha_{w \lambda}$ have alternating signs, in
  other words $(-1)^{|\lambda| - \ell(w)} \alpha_{w \lambda}$ is
  non-negative.
\end{conj}

After the first version of this paper was circulated, we were informed
by A.~Lascoux \cite{lascoux:private} of a beautiful recursive formula
for the coefficients $\alpha_{w \lambda}$ which in particular implies
that \refconj{conj_alpha} is true.

It would be interesting to define a non-commutative polynomial
$G_w(u)$ in the local plactic algebra for any permutation $w$.  Such a
definition might lead to a Littlewood-Richardson rule for the
coefficients $\alpha_{w \lambda}$.

% Since the term of lowest degree in $G_w(x)$ is the stable Schubert
% polynomial or Stanley symmetric function $F_w(x)$, it follows that 
% when $|\lambda| = \ell(w)$ the coefficient $\alpha_{w \lambda}$ is the
% one defined by Stanley \cite{stanley:on*1}.  The coefficients
% $\alpha_{w \lambda}$ also generalize the structure constants
% $c^\nu_{\lambda \mu}$ and $d^\nu_{\lambda \mu}$ of $\Gamma$.

% In the first version of this paper we conjectured that the
% coefficients $\alpha_{w \lambda}$ have alternating signs, i.e.\ 
% $(-1)^{|\lambda| - \ell(w)} \alpha_{w \lambda} \geq 0$.  However,
% A.~Lascoux has sent us a beautiful inductive formula for these
% coefficients which implies this property.  Lascoux's formula is based
% on a general transition formula for Grothendieck polynomials which has
% not been published, but which is related to his results in
% \cite{lascoux:anneau}.

% The coefficients $\alpha_{w \lambda}$ are special cases of quiver
% coefficients $c_r(\mu)$ which will be introduced in
% \cite{buch:structure}.  We also believe that these quiver coefficients
% have alternating signs.

%%% Local Variables: 
%%% mode: latex
%%% TeX-master: "gamma"
%%% End: 

\section{Consequences of the Littlewood-Richardson rule}
\label{sec_conseq}

In this section we will derive some consequences of the formulas
proved in the past two sections.  We will start with a Pieri rule for
the coproduct in $\Gamma$ which is analogous to Lenart's Pieri rule of
\refcor{cor_pieri}.  If $\lambda = (\lambda_1, \dots, \lambda_p)$ is a
partition, we let $\bar \lambda = (\lambda_2, \dots, \lambda_p)$ be
the partition obtained by removing the top row of $\lambda$.

\begin{cor}
\label{cor_copieri}
If $\lambda/\mu$ is a horizontal strip and $k \geq 0$ an integer, then
\[ d^\lambda_{\mu,(k)} = (-1)^{k - |\lambda/\mu|} 
   \binom{r(\mu/\bar \lambda)}{k - |\lambda/\mu|} 
\]
where $r(\mu/\bar \lambda)$ is the number of rows in $\mu/\bar
\lambda$.  If $\lambda/\mu$ is not a horizontal strip then
$d^\lambda_{\mu,(k)} = 0$.
\end{cor}

This is an immediate consequence of \refthm{thm_lrsplit}.  Notice that
this implies that
\[ G_\lambda(x_1,\dots,x_p) = \sum_{\mu,k} (-1)^{k - |\lambda/\mu|}
   \binom{r(\mu/\bar \lambda)}{k - |\lambda/\mu|} \,
   G_\mu(x_1,\dots,x_{p-1}) \, x_p^k
\]
where the sum is over all integers $k \geq 0$ and all partitions $\mu
\subset \lambda$ such that $\lambda/\mu$ is a horizontal strip.  This
gives a practical way to calculate stable Grothendieck polynomials
$G_\lambda(x)$, which can easily be extended to double stable
Grothendieck polynomials and 321-avoiding permutations.  For example,
if $R = (q)^p$ is a rectangle with $p$ rows and $q$ columns, then
\begin{equation}
\label{eqn_rect}
 G_R(x_1,\dots,x_p; y_1,\dots,y_q) = 
   \prod_{\substack{1\leq i\leq p\\ 1\leq j\leq q}} 
   (x_i + y_j - x_i \, y_j) \,.
\end{equation}
To see this, use the Pieri rule and \refcor{cor_factor} to write
\begin{align*} 
  G_R(x_1,\dots,x_p; y) 
  &= \sum_{q \leq j+k \leq q+1} (-1)^{j+k-q} \,
  G_{(q^{p-1},j)}(x_1,\dots,x_{p-1}; y) \cdot x_p^k \\
  &= \sum_{q \leq j+k \leq q+1} (-1)^{j+k-q} \,
  G_{(q)^{p-1}}(x_1,\dots,x_{p-1}; y) \cdot G_j(0; y) \cdot x_p^k \\
  &= G_{(q)^{p-1}}(x_1,\dots,x_{p-1}; y) \cdot G_q(x_p; y)
\end{align*}
where $y$ denotes the variables $y_1, \dots, y_q$.  Since $G_R(x;y) =
G_{R'}(y;x)$, equation (\ref{eqn_rect}) follows from this by induction
on $p$ and $q$.  Notice that $G_1(x_i; y_j) = x_i + y_j - x_i \, y_j$
by \refexm{exm_splitone}, \refthm{thm_321}, and
\reflemma{lemma_conjugate}.

As another application of \refcor{cor_copieri}, notice that if
$\lambda/\mu$ is a horizontal strip then
\[ \sum_{k \geq 0} d^\lambda_{\mu,(k)} =
   \sum_{k \geq 0} (-1)^{k - |\lambda/\mu|} 
     \binom{r(\mu/\bar \lambda)}{k - |\lambda/\mu|} =
   \begin{cases} 1 & \text{if $\mu = \bar \lambda$;} \\
   0 & \text{otherwise.}
   \end{cases}
\]
This implies that
\[ G_\lambda(1, x_1, x_2, \dots) = 
%   \sum_{k \geq 0} G_{\lambda \sslash (k)}(x) \cdot 1^k =
   \sum_{\mu,k} d^\lambda_{\mu,(k)} \, G_\mu(x) \cdot 1^k =
   G_{\bar \lambda}(x) \,.
\]
In particular, if we let $\phi_p : \Gamma \to \Gamma$ be the linear
map which maps $G_\lambda$ to $G_\mu$ where $\mu = (\lambda_{p+1},
\lambda_{p+2}, \dots)$ is obtained by removing the top $p$ rows of
$\lambda$, then $\phi_p$ is a ring homomorphism.  Being used to
calculating in the ring of symmetric functions, we find this somewhat
surprising.  Similarly we have $\Delta \, \phi_p = (\phi_i \otimes
\phi_j) \, \Delta$ whenever $i + j = p$.  This implies lots of
identities among the structure constants $c^\nu_{\lambda \mu}$ and
$d^\nu_{\lambda \mu}$ of $\Gamma$.  For example, if $\emptyset \neq
\mu \subset \lambda$ are partitions and $p \geq \ell(\lambda)$, then
by comparing the coefficients of $G_\mu \otimes 1$ in $\Delta\,
\phi_p\, G_\lambda = (1 \otimes \phi_p)\, \Delta \, G_\lambda$ we
obtain the identity $\sum_\tau d^\lambda_{\mu \tau} = 0$.  The map
$\phi_p$ is related to the proper pushforward map on Grothendieck
groups from a Grassmann bundle to its base variety.  See
\cite{buch:structure} for details about this.

% A similar Pieri rule to \refcor{cor_copieri} is
% $\alpha_{\lambda/(k),\mu} = (-1)^{k-|\lambda/\mu|} \binom{r(\mu/\Bar
% \lambda)-1}{k-|\lambda/\mu|}$ which holds when $\lambda/\mu$ is a
% horizontal strip and can be proved from \refthm{thm_lr321} (this is
% easiest if one starts by conjugating all partitions).

We will next describe some identities involving rectangles.  Given a
rectangle $R = (q)^p$ with $p$ rows and $q$ columns and a Young
diagram $\mu$ contained in $R$, we let $\Hat \mu$ denote $\mu$ rotated
180 degrees and put in the lower-right corner of $R$.  Then
\refthm{thm_lrsplit} implies the following multiplicity free formula
for the coefficients of $\Delta \, G_R$:
\begin{equation}
\label{eqn_splitrect} 
  d^R_{\lambda \mu} = \begin{cases}
  (-1)^{|\lambda| + |\mu| - |R|} & \text{if $\lambda \cup \Hat \mu = R$
  and $\lambda \cap \Hat \mu$ is a rook strip.} \\
  0 & \text{otherwise.}
  \end{cases}
\end{equation}

There is a similar formula for multiplying the stable Grothendieck
polynomials of two rectangles $R_1$ and $R_2$.  Let $R = R_1 \cap R_2$
be their intersection and $\rho = R_1 \cup R_2$ their union.  Then we
have
\begin{equation}
  G_{R_1} \cdot G_{R_2} = \sum_{\lambda,\mu}
  d^R_{\lambda \mu} \, G_{\rho + \lambda, \mu}
\end{equation}
where $(\rho + \lambda, \mu)$ is the partition
\[ (\rho+\lambda,\mu) =
\raisebox{-25pt}{\includegraphics[scale=0.4]{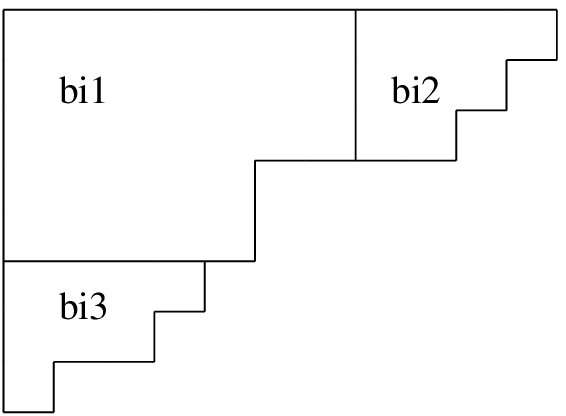}} \]
Assuming that $R_2$ is the tallest of the two rectangles, this
formula is easy to prove from \refthm{thm_lrmult} by counting
set-valued tableaux of shape $R_1 * R_2$.

J.~Stembridge has recently classified all multiplicity free products
of Schur functions \cite{stembridge:multiplicity-free}.  We will say
that a product $G_\lambda \cdot G_\mu$ is multiplicity free if all
constants $c^\nu_{\lambda \mu}$ are $\pm 1$ or zero.  Stembridge's
result then has the following analogue for Grothendieck polynomials.

\begin{prop}
  The product $G_\lambda \cdot G_\mu$ is multiplicity free if and only
  if both partitions $\lambda$ and $\mu$ are rectangles or one of them
  is a single box or empty.
\end{prop}
\begin{proof}
  It is enough to show that if $\lambda$ is not a rectangle and $\ell
  = \ell(\mu) \geq 2$ then the absolute value of some coefficient
  $c^\nu_{\lambda \mu}$ is at least two.  The assumptions on $\lambda$
  and $\ell$ imply that we can find a partition $\rho$ containing
  $\lambda$ such that $\rho/\lambda$ is a vertical strip with $\ell+1$
  boxes and at least three columns.  Then it follows from Lenart's
  Pieri formula that $c^\rho_{\lambda,(1^\ell)} \leq -2$, which means
  that there exist two different tableaux $T_1$ and $T_2$ of shape
  $\lambda$ such that the concatenations $w(T_1) \circ (\ell, \ell-1,
  \dots, 1)$ and $w(T_2) \circ (\ell, \ell-1, \dots, 1)$ are reverse
  lattice words with the same content $\rho$.  But then $w(T_1) \circ
  w(U(\mu))$ and $w(T_2) \circ w(U(\mu))$ must also be reverse lattice
  words with the same content, so $c^\nu_{\lambda,\mu} \leq -2$ if
  $\nu$ is the content of $w(T_i) \circ w(U(\mu))$.
\end{proof}
Similarly, a coproduct $\Delta \, G_\lambda$ is multiplicity free if
and only if $\lambda$ is a rectangle.  

A related consequence of \refthm{thm_lrmult}, which we will need in
\refsec{sec_geometry}, is that if a basis element $G_R$ for a
rectangular partition $R$ occurs with non-zero coefficient in a
product $G_\lambda \cdot G_\mu$, then $R$ is the disjoint union of
$\lambda$ and $\Hat \mu$ (which in particular means that $c^R_{\lambda
  \mu} = 1$).  This can be proved as follows.  Start with $R$ and the
partition $\mu$, and look for a setvalued tableau of shape $\lambda *
\mu$ for some partition $\lambda$, such that this tableau has content
$R$.  The filling of $\mu$ then has to be $U(\mu)$.  Now construct
$\lambda$ and the tableau on $\lambda$ by first filling 1 in some
boxes, then 2, etc.  It is then easy to see that at each step there is
only one choice, i.e.\ both $\lambda$ and the tableau on $\lambda$ are
uniquely determined by the requirement that the word of the tableau on
$\lambda * \mu$ is a reverse lattice word with content $R$.

As noted earlier, the polynomials $G_{\nu/\lambda}$ are uniquely
determined by (\ref{eqn_sslash}).  It is not hard to see that the
formula in the opposite direction is
\begin{equation}
  G_{\nu/\lambda} = \sum_{\sigma \subset \lambda} G_{\nu \sslash \sigma}
\end{equation}
where this sum is over all partitions $\sigma$ contained in $\lambda$.
From this we obtain the following inverse of the relation
$d^\nu_{\lambda \mu} = c^{R+\lambda,\mu}_{\nu,R}$ between the
structure constants of $\Gamma$.  Namely, if $R$ is a rectangle which
contains $\lambda$ and $\mu$ then
\[ c^\nu_{\lambda \mu} = \alpha_{\mu*\lambda,\nu} =
   \alpha_{(R+\lambda,\mu)/R, \nu} =
   \sum_{\sigma \subset R} d^{R+\lambda,\mu}_{\nu,\sigma} \,. 
\]

We will finish this section by proving some results concerning the
shapes of partitions which give non-zero constants $c^\nu_{\lambda
  \mu}$, $d^\nu_{\lambda \mu}$, or $\alpha_{\nu/\lambda,\mu}$.  We
will need a few lemmas which allow us to make small changes to
set-valued tableaux.  Let the {\em shared boxes\/} of a tableau be the
boxes that contain two or more integers.

\begin{lemma}
\label{lemma_skiplarge}
Let $T$ be a set-valued tableau of shape $\nu/\lambda$ with at least
one shared box, and let $y$ be the largest integer contained in a
shared box of $T$.  Suppose $v$ is a sequence of integers such that
the concatenation $w(T) \circ v$ of the word of $T$ with $v$ is a
reverse lattice word.  Then there exists a set-valued tableau $\Tilde
T$ of shape $\nu/\lambda$ such that $w(\Tilde T) \circ v$ is a reverse
lattice word, and so that the integers contained in $\Tilde T$ are the
same as those in $T$, except one integer $x \geq y$ is left out.  (In
other words, the content of $(x) \circ w(\Tilde T)$ is equal to the
content of $w(T)$.)
\end{lemma}
\begin{proof}
  Start by locating the leftmost shared box $A$ of $T$ which contains
  $y$.  To construct $\Tilde T$ we start by removing $y$ from this
  box.  Then look for the nearest box $B$ below or to the left of $A$
  which contains $y+1$, such that the box above $B$ does not contain
  $y$ and the box to the left of $B$ does not contain $y+1$.  If no
  such box exists, then $w(T) \circ v$ stays a reverse lattice word
  even if $y$ is removed from $A$.  If we can locate a box $B$
  satisfying these requirements, we replace $y+1$ with $y$ in this
  box.  Notice that $B$ can't be a shared box by the assumptions.  We
  then continue in the same way, with $B$ in the role of $A$ and $y+1$
  in the role of $y$.  $\Tilde T$ is the tableau resulting when no new
  box $B$ can be obtained.
\end{proof}
The following picture shows an example of the transformation described
in the proof.  The initial box $A$ is the one in the upper-right
corner and $y=6$.
\[ \rtab{-35pt}{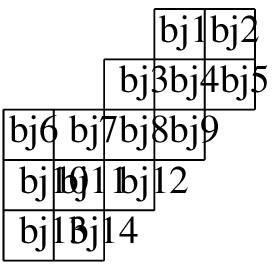} \rightsquigarrow~~~ \rtab{-35pt}{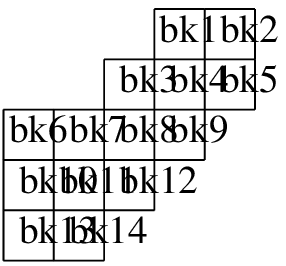} \]

\begin{cor}[of proof]
\label{cor_addbox}
  With the assumptions of \reflemma{lemma_skiplarge}, there exists a
  partition $\rho$ obtained by adding a single box to $\nu$ and a
  tableau $\Tilde T_1$ of shape $\rho/\lambda$ such that $w(\Tilde
  T_1) \circ v$ is a reverse lattice word and the content of $w(\Tilde
  T_1)$ is equal to the content of $w(T)$.
\end{cor}
\begin{proof}
  If $\lambda$ is empty so the shape of $T$ is the partition $\nu$,
  then we obtain $\Tilde T_1$ as the product $x \cdot \Tilde T$ where
  $x$ is the integer which $\Tilde T$ lacks compared to $T$.  When
  this product is formed, the only boxes that can be effected are
  those containing integers strictly larger than $y$ or those
  containing $y$ which are to the left of the original box $A$ in the
  construction of $\Tilde T$.  Since this implies that no shared boxes
  are modified, the shape of $\Tilde T_1$ is only one box larger than
  $\nu$.  When $\nu/\lambda$ is not a partition, the same trick will
  work if we pretend that the boxes of $\lambda$ are actually filled
  with small integers when the product $x \cdot \Tilde T$ is formed.
\end{proof}

Notice that a sequence $w$ of integers between $1$ and $b$ is a
partial reverse lattice word with respect to two integer intervals
$[1,a]$ and $[a+1,b]$ if and only if $w \circ (a^N, \dots, 2^N, 1^N)$
is a reverse lattice word for large $N$.  Therefore
\reflemma{lemma_skiplarge} and \refcor{cor_addbox} are still true if
one replaces ``reverse lattice word'' with ``partial reverse lattice
word'' for given integer intervals.

The following result says that the non-zero structure constants of
$\Gamma$ come in paths which start in the usual Littlewood-Richardson
coefficients.

\begin{prop}
\label{prop_interval}
Let $\lambda$, $\mu$, and $\nu$ be partitions.
\begin{romenum}
\item If $c^\nu_{\lambda \mu} \neq 0$ and $|\nu| > |\lambda| + |\mu|$
  then there exists a partition $\Tilde \nu \subset \nu$ of weight
  $|\Tilde \nu| = |\nu|-1$ such that $c^{\Tilde \nu}_{\lambda \mu}
  \neq 0$.
\item If $c^\nu_{\lambda \mu} \neq 0$ and $|\nu| > |\lambda| + |\mu|$
  then there exists a partition $\Tilde \mu \supset \mu$ of weight
  $|\Tilde \mu| = |\mu|+1$ such that $c^{\nu}_{\lambda \Tilde \mu}
  \neq 0$.
\item If $d^\nu_{\lambda \mu} \neq 0$ and $|\nu| < |\lambda| + |\mu|$
  then there exists a partition $\Tilde \nu \supset \nu$ of weight
  $|\Tilde \nu| = |\nu|+1$ such that $d^{\Tilde \nu}_{\lambda \mu}
  \neq 0$.
\item If $d^\nu_{\lambda \mu} \neq 0$ and $|\nu| < |\lambda| + |\mu|$
  then there exists a partition $\Tilde \mu \subset \mu$ of weight
  $|\Tilde \mu| = |\mu|-1$ such that $d^{\nu}_{\lambda \Tilde \mu}
  \neq 0$.
\item If $\alpha_{\nu/\lambda,\mu} \neq 0$ and $|\mu| > |\nu/\lambda|$
  then there exists a partition $\Tilde \mu \subset \mu$ of weight
  $|\Tilde \mu| = |\mu|-1$ such that $\alpha_{\nu/\lambda,\Tilde \mu}
  \neq 0$.
\item If $\alpha_{\nu/\lambda,\mu} \neq 0$ and $|\mu| > |\nu/\lambda|$
  then there exists a partition $\Tilde \nu \supset \nu$ of weight
  $|\Tilde \nu| = |\nu|+1$ such that $\alpha_{\Tilde \nu/\lambda,\mu}
  \neq 0$.
\item If $\alpha_{\nu/\lambda,\mu} \neq 0$ and $|\mu| > |\nu/\lambda|$
  then there exists a partition $\Tilde \lambda \subset \lambda$ of
  weight $|\Tilde \lambda| = |\lambda|-1$ such that
  $\alpha_{\nu/\Tilde \lambda,\mu} \neq 0$.
\end{romenum}
\end{prop}
\begin{proof}
  \reflemma{lemma_skiplarge} implies (i), (iv), and (v), while
  \refcor{cor_addbox} implies (ii), (iii), and (vi).  For example, to
  prove (ii) from the corollary, recall that if $c^\nu_{\lambda \mu}
  \neq 0$ then there exists a tableau $T$ of shape $\mu$ such that
  $w(T) \circ w(U(\lambda))$ is a reverse lattice word with content
  $\nu$.  Since $|\nu| > |\lambda| + |\mu|$, $T$ must contain a shared
  box.  We can therefore let $\Tilde \mu$ be the shape of the tableau
  $\Tilde T_1$ of \refcor{cor_addbox}.
  
  Finally, to prove (vii) we need to show that if $T$ is a tableau of
  shape $\nu/\lambda$ with at least one shared box such that $w(T)$ is
  a reverse lattice word with content $\nu$, then there exists a
  tableau $\Tilde T$ of a shape $\nu/\Tilde \lambda$ such that
  $w(\Tilde T)$ is a reverse lattice word with the same content $\nu$.
  Let $y$ be the smallest integer contained in a shared box of $T$,
  and let $A$ be the northernmost shared box containing $y$.  Then
  start by removing $y$ from this box.  If all integers in the row
  above $A$ are larger than or equal to $y$, then we can add a new box
  containing $y$ at the left end of this row.  Otherwise let $B$ be
  the rightmost box in the row above $A$ which contains an integer
  strictly less than $y$.  Now replace the integer in $B$ with $y$ and
  continue in the same way with this integer in the role of $y$ and
  $B$ in the role of $A$.  Using the induction hypothesis that some
  box strictly north of $A$ contains $y$, it is not hard to check that
  this process stops before we reach the top row of $T$, and that the
  result is a tableau $\Tilde T$ with the desired properties.
\end{proof}

\begin{lemma}
\label{lemma_skipsmall}
Let $T$ be a set-valued tableau whose shape is a partition $\lambda$,
such that $T$ has at least one shared box.  Let $y$ be the smallest
integer contained in a shared box of $T$.  Suppose $v$ is a sequence
of integers such that $w(T) \circ v$ is a reverse lattice word.  Then
there exists a tableau $\Tilde T$ of shape $\lambda$ and an integer $x
\leq y$ such that $w(\Tilde T) \circ v$ is a reverse lattice word and
the content of $(x) \circ w(\Tilde T)$ is equal to the content of
$w(T)$.
\end{lemma}
\begin{proof}
Let $a$ be the set in the leftmost shared box which contains $y$ and
let $T_1$, $C$, $D$, and $T_2$ be as in the picture.
\[ T = \raisebox{-29pt}{\includegraphics[scale=0.7]{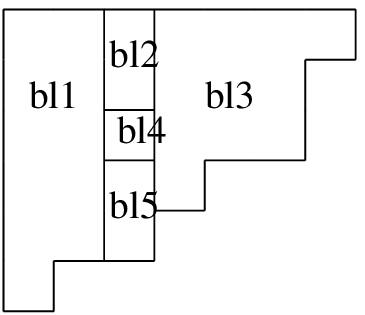}} \]
Let $(T_1, y)$ be the tableau obtained by attaching a box containing
$y$ to the right side of $T_1$ and let $\theta = \sh(T_1,y)/\sh(T_1)$
be the skew diagram of this box.  Then set $(x,\Tilde T_1) =
\RR_\theta(T_1,y)$, and let $\Tilde T$ be the tableau obtained from
$T$ by replacing $T_1$ with $\Tilde T_1$ and $a$ with $\Tilde a = a
\smallsetminus \{y\}$.  Notice that $x$ must be a single integer,
since only integers less than or equal to $y$ in $T_1$ are affected
when forming $\RR_\theta(T_1,y)$, and none of these are in shared
boxes.

Since $w(T) \circ v = w(T_1) \circ w(C) \circ w(a) \circ w(D) \circ
w(T_2) \circ v$ is a reverse lattice word, so is $w(T_1, y) \circ w(C)
\circ w(\Tilde a) \circ w(D) \circ w(T_2) \circ v$.  But then
\reflemma{lemma_revlat} implies that $(x) \circ w(\Tilde T_1) \circ
w(C) \circ w(\Tilde a) \circ w(D) \circ w(T_2) \circ v = (x) \circ
w(\Tilde T) \circ v$ is a reverse lattice word.  The lemma follows
from this.
\end{proof}

\begin{prop}
  If $c^\nu_{\lambda \mu} \neq 0$ then $\nu$ is contained in the union
  of all partitions $\rho$ of weight $|\rho| = |\lambda| + |\mu|$,
  such that the Littlewood-Richardson coefficient $c^\rho_{\lambda
    \mu}$ is non-zero.
\end{prop}
\begin{proof}
  It is enough to show that for each $1 \leq i \leq \ell(\nu)$ there
  is a partition $\rho$ of weight $|\rho| = |\lambda| + |\mu|$ such
  that $c^\rho_{\lambda \mu} \neq 0$ and $\rho_i = \nu_i$.  We will do
  this by induction on $|\nu|$, the case $|\nu| = |\lambda| + |\mu|$
  being trivial.
  
  By \refthm{thm_lrmult} there exists a tableau $T$ of shape $\lambda$
  such that $w(T) \circ w(U(\mu))$ is a reverse lattice word with
  content $\nu$.  Since $|\nu| > |\lambda| + |\mu|$, this tableau must
  have at least one shared box.  If this box contains an integer which
  is larger than $i$, then \reflemma{lemma_skiplarge} gives us a
  tableau $\Tilde T$ in which the number of $i$'s is the same as in
  $T$.  Otherwise some shared box contains an integer smaller than
  $i$, in which case we use \reflemma{lemma_skipsmall} to produce
  $\Tilde T$.  Now if $\Tilde \nu = \sh(\Tilde T)$, then $c^{\Tilde
    \nu}_{\lambda \mu} \neq 0$ and $\Tilde \nu_i = \nu_i$.  Since
  $|\Tilde \nu| = |\nu| - 1$, the required partition $\rho$ exists by
  induction.
\end{proof}

Let us finally remark that the triples of partitions
$(\lambda,\mu,\nu)$ for which $c^\nu_{\lambda \mu} \neq 0$ do not form
a semigroup as is the case if one only considers Littlewood-Richardson
coefficients \cite{zelevinsky:littlewood-richardson}.  For example,
$c^{(2,1)}_{(1),(1)} = -1$ but $c^{(4,2)}_{(2),(2)} = 0$.  However,
$\Gamma$ might still have the property that if $c^{N \nu}_{N \lambda,
  N \mu}$ is not zero for some integer $N > 1$ then $c^\nu_{\lambda
  \mu} \neq 0$.  For Littlewood-Richardson coefficients this has been
proved by Knutson and Tao \cite{knutson.tao:honeycomb}.  The same
question applies to the coefficients $d^\nu_{\lambda \mu}$ and
$\alpha_{\nu/\lambda,\mu}$ as well.

%%% Local Variables: 
%%% mode: latex
%%% TeX-master: "gamma"
%%% End: 

\section{$K$-theory of Grassmannians}
\label{sec_geometry}

This section establishes the link between $K$-theory of Grassmann
varieties and the bialgebra $\Gamma$.  This is then used to describe
some results of A.~Knutson regarding $K$-theoretic triple
intersections on Grassmannians.

% In this last section we will give a discussion of the overall
% structure of $\Gamma$, including its relations to geometry and to the
% ring of symmetric functions.  We will start by making precise what it
% says about $K$-theory of Grassmannians.

If $E$ and $F$ are vector bundles over a variety $X$ and $w$ is a
permutation, we define an element $G_w(F - E)$ in the Grothendieck
group of $X$ as follows.  Suppose first that $E = L_1 \oplus \cdots
\oplus L_e$ and $F = M_1 \oplus \cdots \oplus M_f$ are direct sums of
line bundles.  Then we set
\[ G_w(F - E) = G_w(1-M_1^{-1}, \dots, 1-M_f^{-1} ; 
   1-L_1, \dots, 1-L_e) \,.
\]
Since the stable Grothendieck polynomial $G_w(x;y)$ is symmetric in
both the $x_i$ and the $y_i$ separately, this expression can be
written as a polynomial in the exterior powers of $F^\vee$ and $E$.
For this reason the definition makes sense even when $E$ and $F$ do
not have decompositions into line bundles.  The fact that
$G_\lambda(1-e^{-x}; 1-e^y)$ is super symmetric translates into the
formula $G_\lambda(F \oplus H - E \oplus H) = G_\lambda(F - E)$ where
$H$ is an arbitrary vector bundle.  \reflemma{lemma_conjugate} says
that $G_\lambda(F - E) = G_{\lambda'}(E^\vee - F^\vee)$.

% This notation is motivated by the following $K$-theory version of the
% Thom-Porteous formula (see e.g.\ \cite[Thm.~14.4]{fulton:intersection}
% for the analogue in cohomology).  If $E \to F$ is a general map of
% vector bundles of ranks $e$ and $f$ over a Cohen-Macaulay variety $X$,
% and $\Omega_r(E \to F) \subset X$ is the subscheme of points $x$ for
% which the linear map on fibers $E(x) \to F(x)$ has rank at most $r$,
% then the class of the structure sheaf is given by
% \[ [\O_{\Omega_\lambda(E \to F)}] = G_\lambda(F - E) \]
% where $\lambda = (e-r)^{f-r}$ is the rectangular partition with $f-r$
% rows and $e-r$ columns.  This is a special case of the formula for the
% structure sheaf of a quiver variety which will be proved in
% \cite{buch:structure}.

Now let $X = \Gr(d, \C^n)$ be the Grassmann variety of $d$ dimensional
subspaces of $\C^n$, and let $S \subset \C^n \times X$ be the
tautological subbundle of rank $d$ on $X$.  Let $\lambda$ be a
partition with at most $d$ rows and at most $n-d$ columns.  Then the
class of the structure sheaf of the Schubert variety $\Omega_\lambda
\subset X$ defined by (\ref{eqn_grschub}) is given by
\begin{equation}
\label{eqn_grclass}
  [\O_{\Omega_\lambda}] = G_\lambda(S^\vee) = G_\lambda(S^\vee - 0) \,.
\end{equation}
To see this, let $Y = \Fl(\C^n) = \{ V_1 \subset \cdots \subset V_n =
\C^n \}$ be the variety of full flags in $\C^n$ with tautological flag
$F_1 \subset \dots \subset F_n = \C^n \times Y$.  For any permutation
$w \in S_n$ there is a Schubert variety in $Y$ defined by
\[ \Omega_w = \{ V_\bull \in Y \mid \dim(V_p \cap \C^k) \geq 
   p - r_w(p,n-k) ~\forall p,k \}
\]
where $r_w(p,k) = \# \{ i \leq p \mid w(i) \leq k \}$.  It follows
from \cite{lascoux.schutzenberger:structure} or
\cite[Thm.~3]{fulton.lascoux:pieri} that $[\O_{\Omega_w}] =
\Groth_w(1-L_1, \dots, 1-L_n)$ in $K^\circ Y$, where $\Groth_w(x)$ is
the single Grothendieck polynomial associated to $w$ and $L_i =
F_i/F_{i-1}$.  Now if $\rho : Y \to X$ is the map sending a flag
$V_\bull$ to the subspace $V_d$ of dimension $d$, one can check
\cite[Prop.~10.9]{fulton:young} that $\rho^{-1}(\Omega_\lambda) =
\Omega_{w_\lambda}$, where $w_\lambda \in S_n$ is the Grassmannian
permutation for $\lambda$ with descent in position $d$.  Since the
pullback map $\rho^* : K^\circ X \to K^\circ Y$ is injective, it is
therefore enough to show that $[\O_{\Omega_{w_\lambda}}] =
G_{w_\lambda}(F_d^\vee)$ in $K^\circ(Y)$.  This is true because
$\Groth_{w_\lambda}(x) = G_\lambda(x_1,\dots,x_d)$ by
\refthm{thm_gtos}.

%\[ [\O_{\Omega_{w_\lambda}}] = \Groth_{w_\lambda}(1-L_1, \dots, 1-L_n) =
%   G_\lambda(1-L_1,\dots,1-L_d) = G_\lambda(F_d^\vee) \,. 
%\]

Now given any partition $\nu$, let $I_\nu \subset \Gamma$ be the ideal
spanned by the elements $G_\lambda$ for all partitions $\lambda$ which
are not contained in $\nu$.

\begin{thm}
\label{thm_kgr}
  The map $G_\lambda \mapsto G_\lambda(S^\vee)$ induces an isomorphism
  of rings $\Gamma/I_R \cong K^\circ \Gr(d,\C^n)$ where $R = (n-d)^d$
  is a rectangle with $d$ rows and $n-d$ columns.
\end{thm}
\begin{proof}
  Since the map $G_\lambda \mapsto G_\lambda(S^\vee)$ is surjective by
  \refeqn{eqn_grclass} and since $\Gamma/I_R$ and $K^\circ X$ are free
  Abelian groups of the same rank, it is enough to show that
  $G_\lambda(S^\vee) = 0$ when $\lambda \not \subset R$.  If
  $\ell(\lambda) > d$ this follows from \refthm{thm_321} since $S$ has
  rank $d$.  Now if $0 \to S \to \C^n \to Q \to 0$ denotes the
  universal exact sequence on $X$, we get $G_\lambda(S^\vee) =
  G_\lambda(S^\vee - \C^n) = G_{\lambda'}(\C^n - S) = G_{\lambda'}(Q
  \oplus S - S) = G_{\lambda'}(Q)$ which is zero if $\lambda_1 > n-d =
  \rank(Q)$.
\end{proof}

As mentioned in the introduction, the coproduct on $\Gamma$ is also
closely related to $K$-theory of Grassmannians.  Given positive
integers $d_1 < n_1$ and $d_2 < n_2$, set $X_1 = \Gr(d_1, \C^{n_1})$,
$X_2 = \Gr(d_2, \C^{n_2})$, and $X = \Gr(d_1+d_2, \C^{n_1+n_2})$, and
let $S_1$, $S_2$, and $S$ be the tautological subbundles on these
varieties.  Let $P$ be the product $P = X_1 \times X_2$ with
projections $\pi_i : P \to X_i$, and let $\phi : P \to X$ be the
embedding which maps a pair $(V_1, V_2)$ of subspaces $V_1 \subset
\C^{n_1}$ and $V_2 \subset \C^{n_2}$ to the subspace $V_1 \oplus V_2
\subset \C^{n_1} \oplus \C^{n_2}$.  Then since $\phi^* S = \pi_1^* S_1
\oplus \pi_2^* S_2$ it follows that the pullback on Grothendieck rings
$\phi^* : K^\circ X \to K^\circ P = K^\circ X_1 \otimes K^\circ X_2$
is given by $\phi^*(G_\nu(S)) = G_\nu(\pi_1^* S_1 \oplus \pi_2^* S_2)
= \sum d^\nu_{\lambda \mu} \, G_\lambda(S_1) \otimes G_\mu(S_2)$.

% Triple intersections:

We will next report on some unpublished results of A.~Knutson
regarding triple intersections of Schubert structure
sheaves\footnote{While Knutson's results hold for arbitrary partial
  flag varieties, we shall only be concerned with Grassmannians
  here.}.  Let $\rho : X = \Gr(d, \C^n) \to \{ * \}$ be a map to a
point and let $\rho_* : K^\circ X \to \Z$ be the induced map on
Grothendieck groups.  The triple intersection number of the structure
sheaves $\O_{\Omega_\lambda}$, $\O_{\Omega_\mu}$, and
$\O_{\Omega_\nu}$ is the integer $\rho_*([\O_{\Omega_\lambda}] \cdot
[\O_{\Omega_\mu}] \cdot [\O_{\Omega_\nu}])$.  This is a natural
$K$-theory parallel of the symmetric Littlewood-Richardson
coefficients studied in e.g.\ \cite{knutson.tao:honeycomb}.  Let
$\I_\lambda = \I_{\Omega_\lambda \smallsetminus \Omega^\circ_\lambda}
\subset \O_{\Omega_\lambda}$ denote the ideal sheaf of the complement
of the open Schubert cell $\Omega^\circ_\lambda$ in $\Omega_\lambda$.
To analyze triple intersections, Knutson proved that these ideal
sheaves form a dual basis to the basis of Schubert structure sheaves
with respect to the pairing $(\alpha, \beta) = \rho_*(\alpha \cdot
\beta)$ on $K^\circ X$.  Knutson furthermore worked out the change of
basis matrices.  We will apply the methods developed in the present
paper to give simple proofs of these results.  In addition we will
prove an explicit formula for triple intersections and give an example
showing that these numbers can be negative, although the signs of
triple intersections do not alternate in a simple way.

Let $t = 1 - G_1 \in \Gamma$.  We will abuse notation and write $t$
also for its image $1 - [\O_{\Omega_1}]$ in $K^\circ X$, which by
definition of the polynomial $G_1$ is equal to the class of the line
bundle $\bigwedge^d S$.  \refcor{cor_pieri} implies that for any
partition $\lambda$ we have $t \cdot G_\lambda = \sum
(-1)^{|\sigma/\lambda|} G_\sigma$ where the sum is over all partitions
$\sigma \supset \lambda$ such that $\sigma/\lambda$ is a rook strip.
Since $\rho_*([\O_{\Omega_\sigma}]) = 1$ for each $\sigma \subset R =
(n-d)^d$, it follows from this that $\rho_*(t \cdot
[\O_{\Omega_\lambda}])$ is equal to one when $\lambda = R$ and zero
otherwise.  If $\lambda$ is contained in $R$, let $\Tilde \lambda$ be
the partition obtained by rotating the skew diagram $R/\lambda$ 180
degrees, i.e.\ $\Tilde \lambda = (n-d-\lambda_d, \dots,
n-d-\lambda_1)$.  As we noted in \refsec{sec_conseq}, the coefficient
of $G_R$ in a product $G_\lambda \cdot G_\mu$ is non-zero if and only
if $\mu = \Tilde \lambda$, and in this case we have $c^R_{\lambda \mu}
= 1$.  It follows from this that $\rho_*(t \cdot [\O_{\Omega_\lambda}]
\cdot [\O_{\Omega_\mu}])$ is equal to one if $\mu = \Tilde \lambda$
and zero otherwise.  We conclude that the elements $t \cdot
[\O_{\Omega_\lambda}]$ form a dual basis to the basis of Schubert
structure sheaves, with $[\O_{\Omega_\lambda}]$ and $t \cdot
[\O_{\Omega_{\Tilde \lambda}}]$ dual to each other.  Since
$\Omega_\lambda \smallsetminus \Omega^\circ_\lambda$ is a zero section
of the line bundle $\bigwedge^d S^\vee$ restricted to $\Omega_\lambda$
\cite{schubert:losung}, we finally deduce that $[\I_\lambda] =
[\bigwedge^d S \otimes \O_{\Omega_\lambda}] = t \cdot
[\O_{\Omega_\lambda}]$.

% images in $K^\circ X$ of the elements $H_\lambda = t
% \cdot G_\lambda$ in $\Gamma$ form a dual basis to the basis of
% Schubert structure sheaves.  The fact that $\I_{\Omega_\lambda
%   \smallsetminus \Omega^\circ_\lambda}$ is the image of $H_\lambda$ in
% $K^\circ X$ follows because $\Omega_\lambda \smallsetminus
% \Omega^\circ_\lambda$ is a zero section of the line bundle
% $\bigwedge^d S^\vee$ restricted to $\Omega_\lambda$.

Now, calculating a triple intersection number
$\rho_*([\O_{\Omega_\lambda}] \cdot [\O_{\Omega_\mu}] \cdot
[\O_{\Omega_\nu}])$ is equivalent to expanding $[\O_{\Omega_\lambda}]
\cdot [\O_{\Omega_\mu}]$ in terms of the dual basis $\{ [\I_\sigma]
\}$ and extracting the coefficient of $[\I_{\Tilde \nu}]$.  This is
the same as the coefficient of $G_{\Tilde \nu}$ when the formal power
series $t^{-1} \cdot G_\lambda \cdot G_\mu$ is written as an infinite
linear combination of the basis elements for $\Gamma$.  Notice that
multiplication by $t^{-1}$ takes any basis element $G_\lambda$ to the
sum of all elements $G_\sigma$ for partitions $\sigma$ containing
$\lambda$; this is the inverse operation to multiplication by $t$.  We
therefore obtain the formula
\begin{equation}
  \rho_*([\O_{\Omega_\lambda}] \cdot [\O_{\Omega_\mu}] \cdot
  [\O_{\Omega_\nu}])
  = \sum_{\sigma \subset \Tilde \nu} c^\sigma_{\lambda \mu} \,.
\end{equation}
Alternatively we have $\rho_*([\O_{\Omega_\lambda}] \cdot
[\O_{\Omega_\mu}] \cdot [\O_{\Omega_\nu}]) = \sum_{\sigma \supset
  \lambda} c^{\Tilde \nu}_{\sigma \mu}$.  It turns out that many of
these intersection numbers are non-negative.  For example, when
$\Tilde \nu$ contains the union of all partitions $\sigma$ with
non-zero coefficient $c^\sigma_{\lambda \mu}$, then the intersection
number is equal to one, which follows from the fact that the map
$\phi_p$ of \refsec{sec_conseq} is a ring homomorphism.  Similarly one
can check that all triple intersections on Grassmannians of dimension
smaller than 20 are non-negative.  However, a direct calculation shows
that the coefficient of $G_{5\,4\,3\,1}$ in $t^{-1} (G_{3\,2\,1})^2$
is $-1$.  In other words negative triple intersections can be found on
$\Gr(4,\C^9)$.

% This coefficient can be found using gcalc, by writing:
% low(g[3,2,1]^2, g[5,4,3,1]);

% Note:  t = det S
%
% X = Gr(d, C^n)
% S \subset C^n tautological subbundle.
%
% W = Omega_1 = { S \cap C^{n-d} \neq 0 }
%   = Omega_{d-1}( S -> C^n/C^{n-d} )
%   = Z( \bigwedge^d S -> C )
%
% This shows that  0 --> \bigwedge^d S -> C --> \O_W --> 0
%
% so  t = 1 - [\O_w] = [\bigwedge^d S]  in  K^\circ X.

% End of Triple intersections

Let us remark here that the signs showing up in the structure
constants of $\Gamma$ are to some extend a matter of choice.  To be
precise, all structure constants of $\Gamma$ with respect to the basis
$\{ (-1)^{|\lambda|} \, G_\lambda \}$ are non-negative.  This
viewpoint is equivalent to working with Fomin and Kirillov's
$\beta$-polynomials, with $\beta = 1$
\cite{fomin.kirillov:grothendieck}.  We have chosen to keep a notation
that leads to signs in order to comply with standard definitions and
to honor the fact that the Schubert structure sheaves on a Grassmann
variety do multiply with alternating signs.

%%% Local Variables: 
%%% mode: latex
%%% TeX-master: "gamma"
%%% End: 

\section{The structure of $\Gamma$}
\label{sec_structure}

In this last section we will give a discussion of the overall
structure of $\Gamma$, including its relation to the ring of symmetric
functions.

Recall that the {\em ring of symmetric functions\/} is the span
$\Lambda = \bigoplus_\lambda \Z \cdot s_\lambda$ of all Schur
functions $s_\lambda = s_\lambda(x)$ \cite{macdonald:symmetric*2,
  fulton:young}.  This is in fact a Hopf algebra
\cite[Ch.~1]{de-concini.procesi:quantum}.  Its structure constants are
the Littlewood-Richardson coefficients, i.e.\ 
\[ s_\lambda \cdot s_\mu = \sum c^\nu_{\lambda \mu} \, s_\nu 
   \hspace{15pt} \text{and} \hspace{15pt}
   \Delta \, s_\nu = \sum c^\nu_{\lambda \mu} \, s_\lambda \otimes s_\mu
\]
where the first sum is over partitions $\nu$ and the second over
partitions $\lambda$ and $\mu$, such that $|\nu| = |\lambda| + |\mu|$
in both cases.  The antipode is given by $S(s_\lambda) =
(-1)^{|\lambda|} \, s_{\lambda'}$.

As noted in the introduction, $\Lambda$ is the associated graded
bialgebra to $\Gamma$ with respect to the filtration $\Gamma_p =
\bigoplus_{|\lambda| \geq p} \Z \cdot G_\lambda$.  This is an
immediate consequence of the fact that $c^\nu_{\lambda \mu}$ and
$d^\nu_{\lambda \mu}$ are both equal to the usual
Littlewood-Richardson coefficient when $|\nu| = |\lambda| + |\mu|$.
Furthermore, if we let $\Hat \Gamma$ and $\Hat \Lambda$ be the
completions of $\Gamma$ and $\Lambda$, consisting of infinite linear
combinations of stable Grothendieck polynomials and Schur functions,
respectively, then $\Hat \Gamma \cong \Hat \Lambda$ as bialgebras.
This is true because if we set the variables $y_i$ to zero, then $\Hat
\Gamma$ and $\Hat \Lambda$ both consist of all symmetric power series
in $\Z\llbracket x_1, x_2, \dots \rrbracket$.

Despite from these facts, $\Gamma$ and $\Lambda$ are not isomorphic as
bialgebras themselves.  In fact, there exists no antipode which makes
$\Gamma$ a Hopf algebra.  Recall that an antipode is a linear map $S :
\Gamma \to \Gamma$, such that $S(1) = 1$ and for each non-empty
partition $\nu$ we have $\sum d^\nu_{\lambda \mu} \, S(G_\lambda)
\cdot G_\mu = 0$, or equivalently
\begin{equation}
\label{eqn_antipode}
  \sum_{\lambda \subset \nu} 
  S(G_\lambda) \, G_{\nu \sslash \lambda} = 0
\end{equation}
where $G_{\nu \sslash \lambda}$ is given by (\ref{eqn_defsslash}).
Taking $\nu = (1)$ we get $S(G_1) \cdot (1 - G_1) + 1 \cdot G_1 = 0$
which implies that $S(t) = t^{-1}$ where $t = 1 - G_1$.  Since $t^{-1}
\in \Hat \Gamma$ is equal to the sum of the elements $G_\lambda$ for
all partitions $\lambda$, this is not an element of $\Gamma$.

% It is interesting to compare this identity with the formula $\sum_{n
% \geq 0} (s_1)^n = \sum_\lambda f^\lambda \, s_\lambda$ which holds
% in the completion of the ring of symmetric functions; $f^\lambda$ is
% the number of standard Young tableaux of shape $\lambda$.

However, $\Gamma$ is not far from being a Hopf algebra.  In fact, if
we let $\Gamma_t$ be the localization generated by $\Gamma$ and
$t^{-1}$, then $\Gamma_t$ is a Hopf algebra.  To see this, notice that
(\ref{eqn_sslash}) implies that $G_{\nu \sslash \nu} = t^m$ where $m$
is the number of inner corners of $\nu$, i.e.\ the number of indices
$i$ such that $\nu_i > \nu_{i+1}$.  By (\ref{eqn_antipode}) we
therefore see that an antipode $S : \Gamma_t \to \Gamma_t$ must
satisfy
\begin{equation}
  S(G_\nu) = - t^{-m} \sum_{\lambda \varsubsetneq \nu} 
  S(G_\lambda) \, G_{\nu \sslash \lambda} \,.
\end{equation}
This equation can be used to define $S : \Gamma_t \to \Gamma_t$.  The
obtained antipode is furthermore a ring homomorphism, since it must
agree with the antipode on the ring of symmetric functions extended to
$\Hat \Lambda = \Hat \Gamma$.

Regarding the structure of $\Gamma$ as an abstract ring, we
conjecture:

\begin{conj}
\label{conj_ab}
(a) Any stable polynomial $G_\lambda$ can be written as a polynomial
in the elements $G_R$ for rectangular partitions $R$ contained in
$\lambda$.

(b) The elements $\{ \dots, G_3,\, G_2,\, G_1,\, G_{(1,1)},\,
G_{(1,1,1)}, \dots \}$ corresponding to partitions with only one row
or one column are algebraically independent.
\end{conj}

Part (a) of this conjecture has been verified for all partitions
$\lambda$ of weight at most 9.  For (b), if we define the degree of a
monomial in the $G_k$ and $G_{(1^\ell)}$ to be the total number of
boxes in the partitions of the factors, then all such monomials of
degree up to 9 are linearly independent.  Furthermore, it is not hard
to prove that for any integer $k \geq 2$, the elements $\{G_k, G_1,
G_{(1,1)}, G_{(1,1,1)}, \dots\}$ are algebraically independent.
Namely, if one uses the lexicographic order on partitions, then the
monomials in these elements all have a different maximal partition
$\lambda$ for which the coefficient of $G_\lambda$ is non-zero.

The conjecture has some interesting consequences for the structure of
$\Gamma$.  If (a) is true, then $\Gamma_t$ is generated by the the
elements $G_k$, $G_{(1^\ell)}$ in addition to $t^{-1}$.  To see this,
notice that if $R = (q)^p$ is a rectangular partition with at least
two rows and two columns, and $\lambda = (q^{p-1}, q-1)$ is the
partition obtained by removing the box in the corner of $R$, then
\[ t \cdot G_R = G_1 \, G_\lambda - G_{\lambda+(1)} - G_{\lambda,(1)} +
  G_{\lambda+(1),(1)} \,.
\]
Using (a) this shows that $t \cdot G_R$ can be written as a polynomial
in $G_{q+1}$, $G_{(1^{p+1})}$, and the elements $G_{\Tilde R}$ for
rectangles $\Tilde R$ which are strictly contained in $R$.  This shows
that $G_R$ is in the ring generated by the elements $G_k$,
$G_{(1^\ell)}$, and $t^{-1}$ by induction on the size of $R$.

However, if (b) is true then the elements $G_k$ and $G_{(1^\ell)}$ do
not generate $\Gamma$ as a ring.  In fact, the identity
\[ t \cdot G_{(2,2)} = 
   G_1 \, G_2 + G_1 \, G_{(1,1)} - G_2 \, G_{(1,1)} - (G_1)^3 
\]
implies that $G_{(2,2)}$ can't be written as a polynomial in these
elements if they are algebraically independent.

Geometrically, the fact that $\Gamma$ is a commutative and
cocommutative bialgebra implies that $\Spec \Gamma$ is an Abelian
semigroup scheme.  The existence of an antipode on $\Gamma_t$ means
that the dense open subset $\Spec \Gamma_t$ is a group scheme.
Furthermore, if \refconj{conj_ab} is true then this open subset looks
like an infinite dimensional affine space with a hyperplane removed.

% The element $t$ has other interesting properties.  By the isomorphism
% of \refthm{thm_kgr}, $t$ corresponds to the top exterior power of the
% tautological subbundle on a Grassmann variety.  We encourage the
% reader to work out the nice multiplicity free formulas for the
% products $t \cdot G_\lambda$ and $t^{-1} \cdot G_\lambda$, the later
% of which lives in $\Hat \Gamma$.

% Furthermore it is not
% hard to prove the formulas
% \begin{align*}
% & t \cdot G_\lambda = \sum_{\nu/\lambda \text{ rook strip}}
%   (-1)^{|\nu/\lambda|} \, G_\nu \hspace{20pt} \text{and} \\
% & t^{-1} \cdot G_\lambda = \sum_{\nu \supset \lambda} G_\nu \in 
%   \Hat \Gamma \,.
% \end{align*}

We will finish this paper by raising some additional questions.  First
of all, several people have asked us when a symmetric power series in
$\Z \llbracket x_1,x_2,\dots, y_1, y_2, \dots \rrbracket$ is an
element in $\Gamma$.  Even when the variables $y_i$ are set to zero,
we do not know the answer to this.

In view of \refconj{conj_ab} it would be very interesting to know the
relations between the elements $G_R$ for rectangular partitions.  We
have also been wondering if $\Gamma$ might be a polynomial ring, i.e.\ 
are there algebraically independent elements $h_1, h_2, \dots$ in
$\Gamma$ such that $\Gamma = \Z[h_1, h_2, \dots]$?  We think this is
not the case but haven't been able to prove it.

% It is not hard to see that the single Grothendieck polynomials
% $\Groth_w(x)$ for all permutations $w$ form a basis for the polynomial
% ring $\Z[x_1, x_2, \dots]$.  Linear independence follows from the fact
% that the lowest term of a Grothendieck polynomial $\Groth_w(x)$ is
% equal to the corresponding Schubert polynomial $\Schub_w(x)$.  For
% each integer $n \geq 1$, let $M_n$ denote the linear span of all
% monomials $x_1^{k_1} \, x_2^{k_2} \cdots x_{n-1}^{k_{n-1}}$ for which
% $k_j \leq n-j$ for each $j$.  Since the isobaric divided difference
% operators $\pi_i$ stabilize $M_n$, and since $M_n$ contains the
% Grothendieck polynomial $\Groth_{w_0}(x)$ for the longest permutation
% in $S_n$, we conclude that the Grothendieck polynomials for
% permutations $w \in S_n$ form a basis for this space.

It is not hard to see that the single Grothendieck polynomials
$\Groth_w(x)$ for all permutations $w$ form a basis for the polynomial
ring $\Z[x_1,x_2,\dots]$.  For example, one can check directly that
the Grothendieck polynomials for permutations in $S_n$ give a basis
for the linear span of all monomials $x_1^{k_1} \, x_2^{k_2} \cdots
x_{n-1}^{k_{n-1}}$ for which $k_j \leq n-j$ for each $j$.
Alternatively one can use a stronger result of Lenart
\cite{lenart:noncommutative} which expresses any single Grothendieck
polynomial $\Groth_w(x)$ as an explicit linear combination of Schubert
polynomials $\Schub_{w'}(x)$ with alternating signs, i.e.\ the sign of
the coefficient of $\Schub_{w'}(x)$ is $(-1)^{\ell(w \, w')}$.  On the
other hand, Lascoux has conjectured that each single Schubert
polynomial is a non-negative linear combination of Grothendieck
polynomials.

Now define Grothendieck structure constants $c^w_{u,v} \in \Z$ by
\[ \Groth_u(x) \cdot \Groth_v(x) = 
   \sum_w c^w_{u,v} \, \Groth_w(x) \,. 
\]
These constants are generalizations of the structure constants for
Schubert polynomials as well as the coefficients $c^\nu_{\lambda \mu}$
discussed in this paper.  If $w_\lambda$, $w_\mu$, and $w_\nu$ are
Grassmannian permutations for $\lambda$, $\mu$, and $\nu$ with
descents at the same position, then $c^\nu_{\lambda \mu} =
c^{w_\nu}_{w_\lambda, w_\mu}$.  Based on our results for Grothendieck
polynomials of Grassmannian permutations given in this paper, as well
as on some computational evidence, we pose:

\begin{conj}
The structure constants for single Grothendieck polynomials have
alternating signs, i.e.\ $(-1)^{\ell(u v w)} c^w_{u,v} \geq 0$.
\end{conj}

%%% Local Variables: 
%%% mode: latex
%%% TeX-master: "gamma"
%%% End: 

%\input{overall}

% \bibliographystyle{amsplain}
% \bibliography{../BibTeX/database,gamma}

\providecommand{\bysame}{\leavevmode\hbox to3em{\hrulefill}\thinspace}

%%% Local Variables: 
%%% mode: latex
%%% TeX-master: "gamma"
%%% End: 

%\newpage
%\input{end}
%\input{list}

\end{document}